\documentclass[12pt]{amsart}

\usepackage{tikz}
\usepackage{array}
\usepackage{caption}
\usetikzlibrary{automata}
\usetikzlibrary{positioning}
\usepackage{tikz-network}

\newtheorem{obj}{Observation}[section]
\newtheorem{remark}{Remark}[section]

\newtheorem{theorem}{Theorem}[section]

\newtheorem{lemma}[theorem]{Lemma}

\DeclareMathOperator{\sol}{sol}
\DeclareMathOperator{\nil}{nil}
\DeclareMathOperator{\cyc}{cyc}
\DeclareMathOperator{\girth}{girth}
\DeclareMathOperator{\diam}{diam}
\DeclareMathOperator{\ABC}{ABC}
\DeclareMathOperator{\GA}{GA}
\DeclareMathOperator{\SCI}{SCI}
\DeclareMathOperator{\Nbd}{Nbd}
\DeclareMathOperator{\gen}{gen}
\DeclareMathOperator{\ir}{ir}

\allowdisplaybreaks

\begin{document}
\title[On a bipartite graph defined on groups]{On a bipartite graph defined on groups}

\author[S. Das, A. Erfanian and R. K. Nath]{Shrabani Das, Ahmad Erfanian  and Rajat Kanti Nath*}

\address{S. Das, Department of Mathematical Sciences, Tezpur University, Napaam-784028, Sonitpur, Assam, India.}

\email{shrabanidas904@gmail.com}

\address{A. Erfanian, Department of Pure Mathematics, Ferdowsi University of Mashhad,
	P.O. Box 1159-91775,
	Mashhad, Iran.}

\email{erfanian@um.ac.ir}

\address{R. K. Nath, Department of Mathematical Sciences, Tezpur University, Napaam-784028, Sonitpur, Assam, India.} 
\email{ rajatkantinath@yahoo.com}

\thanks{*Corresponding author}

\begin{abstract}
Let $G$ be a group and $L(G)$ be the set of all subgroups of $G$. We introduce  a bipartite graph $\mathcal{B}(G)$ on  $G$ whose vertex set is the union of  two sets $G \times G$ and $L(G)$, and two vertices $(a, b) \in G \times G$ and $H \in L(G)$ are adjacent if   $H$ is generated by $a$ and $b$. We establish connections between $\mathcal{B}(G)$ and the generating graph of $G$. We also discuss about various graph parameters such as independence number, domination number,  girth, diameter, matching number, clique number, irredundance number, domatic number and minimum size of a vertex cover of $\mathcal{B}(G)$.
We obtain   relations between $\mathcal{B}(G)$ and certain probabilities associated to finite groups. 
We also obtain expressions for various topological indices of $\mathcal{B}(G)$.
Finally, we realize the structures of $\mathcal{B}(G)$ for the dihedral groups of order $2p$ and $2p^2$  and dicyclic groups of order $4p$ and $4p^2$ (where $p$ is any prime) including certain other small order groups.
\end{abstract}

\thanks{ }

\maketitle

	\section{Introduction}
	Finite groups are being characterized through various graphs defined on it for a long time now. There are quite a few graphs whose vertex set contains elements from a group $G$ and  edges defined by properties of $G$. Some popular graphs defined on groups include the commuting graph (pioneered by Brauer and Fowler in \cite{brauer1955groups}), the non-commuting graph (pioneered by Erd{\"o}s and Neumann  \cite{neumann1976problem}), the generating graph (pioneered by Liebeck and Shalev \cite{LS96}), the power graph (pioneered by  Kelarev and Quinn  \cite{KQ2000}), cyclic/non-cyclic graph (pioneered by Abdollahi and Hassanabadi  \cite{AH2007}), nilpotent/non-nilpotent graph (pioneered by Abdollahi and Zarrin \cite{az2010}), solvable/non-solvable graph (pioneered by Hai-Reuven \cite{Hai-Reuven-2013}), and this list has been growing. A survey on the above mentioned graphs defined on groups can be found in \cite{cameron2021graphs}. 
	
	Let $G$ be a  group and $\mathcal{X}$ be a group property, viz. cyclic, abelian, nilpotent, solvable etc. A graph defined on $G$ is called a $\mathcal{X}$ graph of $G$ if the vertex set is $G$ and two distinct vertices $x$ and $y$ are adjacent if $\langle x, y\rangle$ is a $\mathcal{X}$-group. Thus, commuting graph of $G$ is nothing but `abelian graph' of $G$. Recently, the explicit formulas for the number of spanning trees of commuting graphs associated with some specific groups are obtained in \cite{CMMS-2022}. The complement of $\mathcal{X}$ graph is called non-$\mathcal{X}$ graph of $G$.
	Note that the set $\mathcal{X}(G) := \{x \in G : \langle x, y \rangle \text{ is a $\mathcal{X}$-group for all } y \in G\}$ is the set of all dominant vertices of $\mathcal{X}$ graph of $G$. We have  $\mathcal{X}(G) =$ Cyc$(G)$ (the cyclicizer of $G$), $Z(G)$ (the center of $G$), Nil$(G)$ (the hypercenter of $G$)  and Sol$(G)$ (the solvable radical of $G$) when $\mathcal{X}$ represents ``cyclic", ``abelian", ``nilpotent" and ``solvable" respectively. While studying the connectedness, genus and many other graph parameters of   $\mathcal{X}$ graph and non-$\mathcal{X}$ graph often $\mathcal{X}(G)$ is removed from the vertex set. Recently, two counter examples (see \cite{Das-NN-2024, SN-2024}) to a conjecture of Gutman (see \cite[Conjecture 3.1]{Gutman-2011} and \cite{Gutman-78})  regarding the existence of hyperenergetic graphs are provided through non-commuting graphs of finite groups. 
	The power graph of $G$ has vertex set $G$ and two vertices $x$ and $y$
	are adjacent if one of $x$ and $y$ is the power of the other. 
	Recent works on power graph can be found in \cite{KSCC-2021, MPS23, PPS21}.
	The generating graph of a finite group $G$, denoted by $\Gamma_{\gen}(G)$, is a simple undirected graph with vertex set $V(\Gamma_{\gen}(G))= G$ and two distinct vertices $x$ and $y$ are adjacent if $\langle x, y \rangle = G$. There are other classes of graphs defined on groups whose vertices are the orders of the elements (see \cite{MN-2024}) or the orders of the conjugacy classes (see \cite{Lewis-2008}). A survey on graphs whose vertex set consists of conjugacy classes of a group can be found in \cite{CJSN-2024}.
	Another class of graphs have been defined on groups by considering the vertex set as the set of subgroups of the group. For instance, intersection graph (introduced by Cs$\acute{\rm a}$k$\acute{\rm a}$ny and Poll$\acute{\rm a}$k  \cite{CP69}), inclusion graph (introduced by Devi and Rajkumar  \cite{DR16}) and permutability graph (introduced by Rajkumar and Devi  \cite{RD14})  of $G$ are such graphs. 
	In these graphs, if $H$ and $K$ are two vertices (subgroups of $G$) then    $H$ and $K$ are adjacent if (i) $H \cap K \neq \{1\}$ (in case of intersection graph) (ii)  $H \subset K$ or $K \subset H$ (in case of inclusion graph) (iii)  $HK = KH$ (in case of permutability graph).
	Various results on these graphs can be found in \cite{ahm2015,SK-2017,OWW20,RD16,RDG16}.

	Note that none of the above mentioned graphs are bipartite in nature, while the use of bipartite graphs in solving real-world problems has been known.
	Bipartite graphs are extensively used in modern coding theory while receiving code words  from a channel. They can be used in medical field in the detection of lung cancer, throat cancer etc. In economics, one can see how markets work when buyers and sellers do business. Bipartite graphs are also used for optimizing transportation systems, route planning, and logistics management. Reader may have a look to  \cite{Arun-Koma-15} and the references therein for these applications of bipartite graphs.
	
	In this paper, we introduce  a bipartite graph $\mathcal{B}(G)$ on a group $G$ whose vertex set $V(\mathcal{B}(G))$ is the union of two sets $G \times G$ and $L(G) := \{H : H \text{ is a subgroup of } G\}$, and two vertices $(a, b) \in G \times G$ and $H \in L(G)$ are adjacent if  $H = \langle a, b \rangle$, the subgroup generated by $a$ and $b$. We write $V(\mathcal{B}(G)) = G \times G \sqcup L(G)$, where $\times$ denotes the Cartesian product of sets and $\sqcup$ denotes the union of disjoint sets. We shall also use $\sqcup$ to denote the  union of disjoint graphs. The neighborhood of any vertex $x$ of $\mathcal{B}(G)$ is denoted by $\Nbd_{\mathcal{B}(G)}(x):= \{y \in V(\mathcal{B}(G)) : y \text{ is adjacent to } x\}$. Further, if $S$ is a subset of $V(\mathcal{B}(G))$ then we write $\mathcal{B}(G)[S]$ to denote the induced subgraph of $\mathcal{B}(G)$ induced by $S$.
	
	For any   group $G$, it is obvious that every element of $G \times G$ is adjacent to some element of $L(G)$. Also, if $G$ is a $2$-generated group then every element of $L(G)$ is adjacent to some elements of $G \times G$. We also have the following  observations.
	\begin{obj}\label{vrtex_deg_of_X_in_B(G)} Let $G$ be any group. Then $\mathcal{B}(G)$ has the following properties.
		\begin{enumerate}
			\item  For all $(a, b) \in G \times G$, the degree of $(a, b)$ in $\mathcal{B}(G)$  is one. That is, $\deg_{\mathcal{B}(G)}((a, b))$ $ = 1$.

			\item $\mathcal{B}(G)$ has no cycle and it is a forest having $|L(G)|$ components. In particular, $\mathcal{B}(G)$ is a union of \, $|L(G)|$ \, star graphs.    If $G$ is a cyclic group of  prime order then $\mathcal{B}(G) = K_2 \sqcup K_{1, |G|^2-1}$, where $K_n$ and $K_{1, n}$ denote complete graphs on $n$ vertices and star graphs on $n+1$ vertices respectively. 
			
			Let $p$ be any prime and $G = \langle a \rangle$. If  $|G|=2p$ then $V(\mathcal{B}(G)) = G \times G \sqcup \{\{1\}, \langle a^p \rangle, \langle a^2 \rangle, \langle a \rangle\}$. Since \quad   $|\langle a^p \rangle| = 2$ and  $|\langle a^2 \rangle| =p$ we have $\mathcal{B}(G)[\{\langle a^p \rangle\}$ $ \sqcup \Nbd_{\mathcal{B}(G)}(\langle a^p \rangle)] = K_{1, 3}$ and $\mathcal{B}(G)[\{\langle a^2 \rangle\} \sqcup \Nbd_{\mathcal{B}(G)}(\langle a^2 \rangle)] = K_{1, p^2 - 1}$. Also, $\mathcal{B}(G)[\{\langle a \rangle\} \sqcup \Nbd_{\mathcal{B}(G)}(\langle a \rangle)] = K_{1, 3p^2 - 3}$ noting that \quad $|\Nbd_{\mathcal{B}(G)}(\langle a \rangle)| = 4p^2  - $ $(1 + 3 + p^2 - 1) = 3p^2 - 3$. Thus, $\mathcal{B}(G) = K_2 \sqcup K_{1, 3}\sqcup K_{1, p^2 - 1} \sqcup K_{1, 3p^2 - 3}$.
			
			If $|G|=p^2$ then $V(\mathcal{B}(G)) = G \times G \sqcup \{\{1\}, \langle a^p \rangle,  \langle a \rangle\}$. Since  $|\langle a^p \rangle| =p$ we have $\mathcal{B}(G)[\{\langle a^p \rangle\} \sqcup \Nbd_{\mathcal{B}(G)}(\langle a^p \rangle)] = K_{1, p^2 - 1}$. Also, $\mathcal{B}(G)[\{\langle a \rangle\} \sqcup \Nbd_{\mathcal{B}(G)}(\langle a \rangle)] = K_{1, p^4 - p^2}$ noting that  $|\Nbd_{\mathcal{B}(G)}(\langle a \rangle)| = p^4 - (1 + p^2 - 1)$. Thus, $\mathcal{B}(G) = K_2 $ $\sqcup  K_{1, p^2 - 1} \sqcup K_{1, p^4 - p^2}$.
			
			If $|G|= 2p^2$ (for odd prime $p$) then $V(\mathcal{B}(G)) = G \times G$ $ \sqcup \{\{1\}, \langle a^{p^2} \rangle, \langle a^{2p} \rangle, \langle a^p \rangle$, $\langle a^2 \rangle, \langle a \rangle\}$. 
			Since $|\langle a^{p^2} \rangle| = 2$,  $|\langle a^{2p} \rangle| = p$, $|\langle a^p \rangle|= 2p$ and  $|\langle a^2 \rangle| =p^2$ we have $\mathcal{B}(G)[\{\langle a^{p^2} \rangle\} \sqcup \Nbd_{\mathcal{B}(G)}(\langle a^{p^2} \rangle)] = K_{1, 3}$,   $\mathcal{B}(G)[\{\langle a^{2p} \rangle\} \sqcup \Nbd_{\mathcal{B}(G)}(\langle a^{2p} \rangle)] = K_{1, p^2-1}$, 
			$\mathcal{B}(G)[\{\langle a^p \rangle\} \sqcup \Nbd_{\mathcal{B}(G)}(\langle a^p \rangle)] = K_{1, 3p^2 - 3}$, $\mathcal{B}(G)[\{\langle a^2 \rangle\} \sqcup \Nbd_{\mathcal{B}(G)}(\langle a^2 \rangle)] = K_{1, p^4 - p^2}$. Also, $\mathcal{B}(G)[\{\langle a \rangle\} \sqcup \Nbd_{\mathcal{B}(G)}(\langle a \rangle)] = K_{1, 3p^4 - 3p^2}$ noting that $|\Nbd_{\mathcal{B}(G)}(\langle a \rangle)| = 4p^4 - (1 + 3 + p^2 - 1 + 3p^2 - 3 + p^4 - p^2)$ $ = 3p^4 - 3p^2$. Thus, $\mathcal{B}(G) = K_2 \sqcup K_{1, 3} \sqcup K_{1, p^2 - 1} \sqcup K_{1, 3p^2 - 3} \sqcup K_{1, p^4 - p^2} \sqcup K_{1, 3p^4 - 3p^2}$.
			
			If $G$ is a non-cyclic group of order $p^2$ then $G$ has one subgroup of order one, $p + 1$ subgroups of order $p$ and one subgroup of order $p^2$. Let $I = \{1\}$, $H_1, H_2, \dots, H_{p+1}$  and $K = G$ be the subgroups of $G$, where $H_i \cong \mathbb{Z}_p$ for $1 \leq i \leq p+1$. Then $\mathcal{B}(G)[\{I\}\sqcup \Nbd_{\mathcal{B}(G)}(I)] = K_2$, \quad  $\mathcal{B}(G)[\{H_i\}\sqcup \Nbd_{\mathcal{B}(G)}(H_i)] = \mathcal{B}(\mathbb{Z}_p)[\{\mathbb{Z}_p\}\sqcup \Nbd_{\mathcal{B}(\mathbb{Z}_p)}(\mathbb{Z}_p)] = K_{1, p^2 - 1}$ for $1 \leq i \leq p+1$. Further, $\mathcal{B}(G)[\{G\}$ $ \sqcup \Nbd_{\mathcal{B}(G)}(G)]$ $ = K_{1, p(p-1)(p^2 - 1)}$ noting that  $|\Nbd_{\mathcal{B}(G)}(G)| = p^4 - (p+1)(p^2 -1) - 1 = p(p-1)(p^2 - 1)$. Thus, $\mathcal{B}(G) = K_2 \sqcup (p+1)K_{1, p^2 - 1} \sqcup K_{1, p(p-1)(p^2 - 1)}$,  where $mK_{1, n}$ denotes the disjoint union of $m$ copies of the star $K_{1, n}$.
			
			\item     $\mathcal{B}(G)$ is connected if and only if $G = \{1\}$. In this case, $\mathcal{B}(G)= K_2$.
		\end{enumerate}
	\end{obj}

	In Section 2, we obtain some properties of $\mathcal{B}(G)$. In particular, we establish connections between $\mathcal{B}(G)$ and $\Gamma_{\gen}(G)$. We also discuss about various graph parameters such as independence number, domination number,  girth, diameter, matching number, clique number, irredundance number, domatic number and minimum size of a vertex cover of $\mathcal{B}(G)$.
	One big motivation in defining the graph $\mathcal{B}(G)$ is to obtain various probabilities associated to finite groups through this graph. In Section 3, we obtain relations between $\mathcal{B}(G)$ and certain probabilities associated to finite groups. Using those relations, we calculate the exact probabilities for some well-known small order finite groups. We shall also obtain expressions for various topological indices such as first and second Zagreb indices, Randic Connectivity index, Atom-Bond Connectivity index, Geometric-Arithmetic index, Harmonic index and Sum-Connectivity index  of $\mathcal{B}(G)$.
	In Section 4, we first realize the structures of $\mathcal{B}(G)$ when $G = S_3, D_8, Q_8, D_{10}, D_{12}, A_4$ and $S_4$.  After that we   realize the structures of $\mathcal{B}(G)$ when $G = D_{2p}$ and $D_{2p^2}$ the dihedral groups of order $2p$ and $2p^2$ for any prime $p$, where $D_{2n}$ is the dihedral group presented by $\langle a, b: a^n=b^2=1, bab=a^{-1} \rangle$.  We conclude the paper realizing the structures of $\mathcal{B}(G)$ when $G = Q_{4p}$ and $Q_{4p^2}$ the dicyclic groups of order $4p$ and $4p^2$ for any prime $p$, where $Q_{4n}$ is the dicyclic group presented by $\langle a, b : a^{2n} = 1, b^2 = a^n, bab^{-1} = a^{-1} \rangle$.
	
	\section{Some properties of \, $\mathcal{B}(G)$}
	
	We begin with the following properties of $\mathcal{B}(G)$.
	\begin{theorem}
		If $G$ is a non-trivial finite group, then $\deg_{\mathcal{B}(G)}(x) \leq |G|^2-1$ for all $x \in V(\mathcal{B}(G))$. Further $\deg_{\mathcal{B}(G)}(G) = |G|^2-1$ if and only if $G$ is a cyclic group of prime order.
	\end{theorem}
	\begin{proof}
		We have $V(\mathcal{B}(G))=G \times G \, \sqcup \, L(G)$ 
		and $\deg_{\mathcal{B}(G)}(a, b) =1$ for all $(a, b) \in G \times G$. Also, $\{ 1 \} \in L(G)$ and $\{1\}$ is adjacent to $(1, 1)$ only. Therefore, for all $x \in L(G)\setminus \{1\}$, we have $\deg_{\mathcal{B}(G)}(x) \leq |G|^2-1$.
		
		If $G$ is a cyclic group of prime order, then all the non-identity elements of $G$ are its generators. Also, $L(G)=\{\{1\}, G\}$. As such, $\deg_{\mathcal{B}(G)}(\{1\})=1$ since $\{1\}$ is adjacent to $(1, 1)$ only and $\deg_{\mathcal{B}(G)}(G)=|G|^2-1$. Conversely, suppose that $\deg_{\mathcal{B}(G)}(G)=|G|^2-1$. Then  for every element $(1, 1) \ne (a, b) \in G \times G$ we have $\langle a, b\rangle = G$. In particular, $\langle a\rangle = G$ for all $1\ne a \in G$. This shows that $G$ is cyclic group of prime order.                
	\end{proof}
	
	In the following theorem we obtain degree of any vertex   $H \in L(G)$ in the graph $\mathcal{B}(G)$ using the size of the generating graph $\Gamma_{\gen}(H)$. 
	
	\begin{theorem}\label{relatn B(G) and generating graph}
		Let $G$ be a finite group and  $H \in L(G)$.   Then 
		\[
		\deg_{\mathcal{B}(G)}(H)=\begin{cases}
			1, & \text{ if } H=\{1\} \\
			2|e(\Gamma_{\gen}(H))|+\phi(|H|), & \text{ if } H \text{ is cyclic } \\
			2|e(\Gamma_{\gen}(H))|, & \text{ otherwise. }
		\end{cases}
		\]
		Here, $\Gamma_{\gen}(H)$ is the generating graph of $H$ and  $\phi(|H|)$ is the number of generators of $\mathbb{Z}_{|H|}$.
	\end{theorem}
	\begin{proof}
		Clearly, $(1,1)$ is the only vertex adjacent to $\{1\}$ in $\mathcal{B}(G)$ and so $\deg_{\mathcal{B}(G)}(H)=1$ if $H=\{1\}$.
		
		If $H \ne \{1\}$ is a cyclic group then $\phi(|H|)$ gives the number of generators of $H$. We have
		\begin{align*}
			\deg_{\mathcal{B}(G)}(H)&=\left|\{(a,b) \in G \times G: \langle a,b \rangle =H\}\right| \\
			&=\phi(|H|)+\left|\{(a,b) \in G \times G: \langle a,b \rangle =H, a \neq b\}\right|.
		\end{align*}
		Now, for $a \neq b$, if $\langle a,b \rangle=\langle b,a \rangle=H$ then $(a,b)$ and $(b,a)$ are adjacent to $H$ in $\mathcal{B}(G)$ and
		$a$ is adjacent to $b$ in $\Gamma_{\gen}(H)$. It follows that,  the pairs $(a,b), (b,a), a \neq b$ that generates $H$, contribute   one edge in $\Gamma_{\gen}(H)$ and two edges in $\mathcal{B}(G)$. Therefore, $|e(\Gamma_{\gen}(H))|=\frac{1}{2}\left|\{(a,b) \in G \times G: \langle a,b \rangle =H,\right. $ $\left. a \neq b\}\right|$. Thus, $\deg_{\mathcal{B}(G)}(H)=2|e(\Gamma_{\gen}(H))|+\phi(|H|)$.
		
		If $H$ is non-cyclic then \quad	
		$\deg_{\mathcal{B}(G)}(H)=\left|\{(a,b) \in G \times G: \langle a,b \rangle =H, a \neq b\}\right|$, since $\{(a, a) \in G \times G: \langle a, a \rangle =H\}$ is an empty set. 
		Therefore, by  similar arguments as above, it follows that $\deg_{\mathcal{B}(G)}(H)=2|e(\Gamma_{\gen}(H))|$.
	\end{proof}

	The following theorem is useful in obtaining   independence and domination number of $\mathcal{B}(G)$.

	\begin{theorem}\label{size of A bigger than that of B}
		For any $2$-generated finite group $G$, if $A=G \times G$ and $B=L(G)$, then $|A| \geq |B|$ with equality when $G$ is a group of order $1$.
	\end{theorem}
	\begin{proof}
		Define a map $f:G \times G \rightarrow L(G)$ by $f((a, b))= \langle a, b \rangle$ for all $a, b \in G$. We have $f((1, 1))=\{1\}$ and $f((a, 1))=f((1, a))=f((a, a))=\langle a \rangle$ for all $a \in G$. So $f$ is a many-one function. Also, if $G$ is a 2-generated group, then any $H \in L(G)$ is adjacent to some elements of $G \times G$. As such, $f$ is an onto function. Therefore, $|G \times G| > |L(G)|$ when $|G| > 1$. 
		
		For $|G|=1$, $G \times G$ and $L(G)$ have same cardinality equal to one.
	\end{proof}
	
	Let $\Gamma$ be any graph. An independent vertex set of $\Gamma$ is a subset of the vertex set of $\Gamma$ such that no two vertices in the subset represent an edge of $\Gamma$. The cardinality of the largest independent vertex set of $\Gamma$ is the independence number of $\Gamma$. 
	A subset $S$ of $V(\Gamma)$  is said to be a dominating set of vertices in $\Gamma$ if $\left(\cup_{s \in S}\Nbd(s)\right) \cup S = V(\Gamma)$. A dominating set of smallest size is called a minimum dominating set and its cardinality is called the domination number of $\Gamma$.

	\begin{theorem}\label{independence-domination no. of B(G)}
		If $A=G \times G$ and $B=L(G)$, then independence and domination number of $\mathcal{B}(G)$ are the sizes of $A$ and $B$ respectively where $G$ is any $2$-generated finite group.
	\end{theorem}
	\begin{proof}
		By Theorem \ref{size of A bigger than that of B} we have $|A| \geq |B|$. Since $\mathcal{B}(G)$ is a bipartite graph,  by  definition, the independence number of $\mathcal{B}(G)$ is $|A|$. Also, every element of $A$ is adjacent to some elements of $B$ and if $G$ is a $2$-generated finite group, then any element of $B$ is adjacent to some elements of $A$. Therefore, by definition, the domination number of $\mathcal{B}(G)$ is $|B|$.
	\end{proof}
	\begin{remark}\label{first remark for r-generated}
		Let $G$ be an $r$-generated finite group where $r\geq3$.
		\begin{enumerate}
			\item  It can be easily seen that a vertex $H \in B=L(G)$ is isolated in $\mathcal{B}(G)$ if $H$ is generated by 3 or more elements. For example, if $G= \mathbb{Z}_{4} \times  \mathbb{Z}_{4} \times  \mathbb{Z}_{4}$ then $\langle \Bar{2} \rangle \times \langle \Bar{2} \rangle \times \langle \Bar{2} \rangle, \langle \Bar{2} \rangle \times\langle \Bar{2} \rangle \times \langle \mathbb{Z}_{4} \rangle, \langle \Bar{2} \rangle \times \langle \mathbb{Z}_{4} \rangle \times \langle \Bar{2} \rangle, \langle \mathbb{Z}_{4} \rangle \times \langle \Bar{2} \rangle \times \langle \Bar{2} \rangle, \langle \mathbb{Z}_{4} \rangle \times \langle \mathbb{Z}_{4} \rangle \times \langle \Bar{2} \rangle, \langle \mathbb{Z}_{4} \rangle \times \langle \Bar{2} \rangle \times \langle \mathbb{Z}_{4} \rangle, \langle \Bar{2} \rangle \times \langle \mathbb{Z}_{4} \rangle \times \langle \mathbb{Z}_{4} \rangle, \mathbb{Z}_{4} \times  \mathbb{Z}_{4} \times  \mathbb{Z}_{4}$ etc. are some isolated vertices in $\mathcal{B}(G)$. We also have that $|A| =|G \times G| =4096 \geq 129=|L(G)|=|B|$. Thus the conclusion of Theorem \ref{size of A bigger than that of B} is true for $G=\mathbb{Z}_{4} \times  \mathbb{Z}_{4} \times  \mathbb{Z}_{4}$. In general, the conclusion of Theorem \ref{size of A bigger than that of B} may also be true for any finite $r$-generated group where $r\geq 3$. However, the proof we have given will not work in this case as there are isolated vertices.
			\item Let $L_2(G)=\{H \in L(G): H \text{ is generated by 1 or 2 elements} \}$. Then $|A|\geq |L_2(G)|$ and $A \sqcup (L(G) \setminus L_2(G))$ is the largest independent set of $\mathcal{B}(G)$. Hence, independence number of $\mathcal{B}(G)$ is $|A|+|L(G)|-|L_2(G)|$. Further, if $|A|\geq |B|$ then domination number of $\mathcal{B}(G)$ is $|B|$.
		\end{enumerate}
	\end{remark}

	Let $\Gamma$ be any graph. The girth of $\Gamma$, denoted by $\girth(\Gamma)$, is the size of the smallest cycle in it. The diameter of $\Gamma$, denoted by $\diam(\Gamma)$, is defined as the maximum distance of  any two vertices of it. 
	A matching in $\Gamma$ is a subset of the edge set of $\Gamma$ such that no two edges in the subset share common vertices. A maximum matching is a matching that contains the largest possible number of edges. The number of edges in a maximum matching of $\Gamma$ is called the matching number, denoted by $\nu(\Gamma)$, of $\Gamma$. A clique of $\Gamma$ is defined as a subset of $V(\Gamma)$ such that every two distinct vertices of the subset are adjacent. A maximum clique is a clique such that there is no clique with more vertices. The number of vertices in a maximum clique of $\Gamma$ is called the clique number, denoted by $\omega(\Gamma)$, of $\Gamma$. 
	The bondage number of $\Gamma$, denoted by $b(\Gamma)$, is the cardinality of the smallest set $E$ of edges such that the domination number of $\Gamma$ after removing the edges in $E$  is strictly greater than that of original $\Gamma$.

	A subset $S$ of $V(\Gamma)$ is said to be an irredundant set of $\Gamma$ if $\left(\cup_{s \in S \setminus \{v\}}\Nbd(s)\right) \cup \left(S \setminus \{v\}\right) \neq \left(\cup_{s \in S}\Nbd(s)\right) \cup S$, for every vertex $v \in S$. A maximal irredundant set of $\Gamma$ is an irredundant set that cannot be expanded to another irredundant set by addition of any vertex of $\Gamma$. The irredundance number, denoted by $\ir(\Gamma)$, is the minimum size of a maximal irredundant set of $\Gamma$. 
	A domatic partition of $\Gamma$ is a partition of $V(\Gamma)$ into disjoint sets $V_1, V_2, \ldots, V_k$ such that each $V_i$ is a dominating set of $\Gamma$. The maximum size of a domatic partition is called domatic number of $\Gamma$, denoted by $d(\Gamma)$. 
	A vertex cover of $\Gamma$ is a set of vertices of $\Gamma$ that includes at least one endpoint of every edge of $\Gamma$. We write $\beta(\Gamma)$ to denote the minimum size of a vertex cover of $\Gamma$. Note that $\alpha(\Gamma) + \beta(\Gamma) = |V(\Gamma)|$ (see \cite[Corollary 7.1]{BM1977}).
	
	We obtain all the above mentioned graph parameters for  $\mathcal{B}(G)$ in the following result.
	\begin{theorem}
		For any group $G$, the graph $\mathcal{B}(G)$ has the following properties:
		\begin{enumerate}
			\item $\girth(\mathcal{B}(G))= 0$ and $\diam(\mathcal{B}(G)) =1$ or $\infty$.
			\item $\nu(\mathcal{B}(G))=|L_2(G)|$.
			\item $\omega(\mathcal{B}(G))  = 2$.
		\end{enumerate}
	\end{theorem}
	\begin{proof}
		\begin{enumerate}
			\item By Observation \ref{vrtex_deg_of_X_in_B(G)}(b), it follows that $\mathcal{B}(G)$ has no cycle. Therefore, $\girth(\mathcal{B}(G))= 0$. The second part follows from Observation \ref{vrtex_deg_of_X_in_B(G)}(c).

			\item Note that every edge of $\mathcal{B}(G)$ is incident to some $H \in L_2(G)$. Consider a subset $E$ of $e(\mathcal{B}(G))$, the edge set of $\mathcal{B}(G)$, such that there exists only one edge in $E$ with $H$ as its endpoint for each $H \in L_2(G)$. Clearly, $E$ is a matching in $\mathcal{B}(G)$ and $|E|=|L_2(G)|$. Now, if we include any more edge to $E$, there will be two edges in $E$ having a common endpoint $H$ for some $H \in L_2(G)$. Therefore, $E$ is a maximum matching and $\nu(\mathcal{B}(G))=|E|$. Hence, the result follows.
			
			\item Note that $(1, 1)$ and $\{1\}$ are always adjacent in $\mathcal{B}(G)$. Therefore, $\mathcal{B}(G)$ has a clique of size two. Thus, $\omega(\mathcal{B}(G)) \geq 2$. Suppose that $\omega(\mathcal{B}(G)) \geq 3$. Then $\mathcal{B}(G)$ must have a cycle of length greater than or equal to three, which is a contradiction.  Therefore, $\omega(\mathcal{B}(G))=2$.    
		\end{enumerate}
	\end{proof}
	\begin{theorem}
		For a $2$-generated group $G$, we have  $b(\mathcal{B}(G)) =1$,  $\ir(\mathcal{B}(G))=\beta(\mathcal{B}(G))=|L(G)|$ and $d(\mathcal{B}(G))=2$.  
	\end{theorem}
	\begin{proof}
		By Theorem \ref{independence-domination no. of B(G)}, we have domination number of $\mathcal{B}(G)=|B|=|L(G)|$. Also, $\deg_{\mathcal{B}(G)}((a,b))=1$ for any $(a,b)\in G \times G$. If we remove any edge from $\mathcal{B}(G)$, $L(G)$ will not be a dominating set and any other dominating set will have size at least one more than $|L(G)|$. This increases the domination number of the new graph by at least one. Therefore, $b(\mathcal{B}(G))=1$. 
		
		By definition, we have $G \times G$ and $L(G)$ both are maximal irredundant sets of $\mathcal{B}(G)$. From Theorem \ref{size of A bigger than that of B}, we know $|G|^2 \geq |L(G)|$. Therefore, $\ir(\mathcal{B}(G))=|L(G)|$.
		
		We have $\alpha(\mathcal{B}(G))+\beta(\mathcal{B}(G))=|V(\mathcal{B}(G))|$, where $\alpha(\mathcal{B}(G))$ is the independence number of $\mathcal{B}(G)$. From Theorem \ref{independence-domination no. of B(G)}, we have $\alpha(\mathcal{B}(G))=|G|^2$. Therefore $\beta(\mathcal{B}(G))=|G|^2+|L(G)|-|G|^2=|L(G)|$.
		
		We have $V(\mathcal{B}(G))$ is the disjoint union of $G \times G$ and $L(G)$. Also, both $G \times G$ and $L(G)$ are dominating sets of $\mathcal{B}(G)$. As such, $d(\mathcal{B}(G))\geq 2$. It was shown in \cite{CH-1977} that 
		$d(\Gamma) \leq \delta(\Gamma) + 1$ for any graph $\Gamma$, where $\delta(\Gamma)$ is the minimum degree of  $\Gamma$. In our case, $\delta(\mathcal{B}(G)) = 1$ and so $d(\mathcal{B}(G)) \leq 1+1=2$. Hence, $d(\mathcal{B}(G))=2$.
	\end{proof}
	\begin{remark}
		Let $G$ be an $r$-generated group where $r \geq 3$. 
		Then
		\begin{enumerate}
			\item The domination number of $\mathcal{B}(G)$ is \, $\min\left\{|G|^2+|L(G)|-|L_2(G)|, |L(G)|\right\}$. Removing any edge from $\mathcal{B}(G)$ will increase the domination number strictly by one. Therefore, $b(\mathcal{B}(G))=1$.
			\item Both $G \times G \sqcup (L(G) \setminus L_2(G))$ and $L(G)$ are maximal irredundant sets. If $|G|^2 \geq |L(G)|$ then $\ir(\mathcal{B}(G))=|L(G)|$. In general, $$\ir(\mathcal{B}(G))=\min\left\{|G|^2+|L(G)|-|L_2(G)|, |L(G)|\right\}.$$
			\item From Remark \ref{first remark for r-generated}(b), we have independence number of $\mathcal{B}(G)$ is $|G|^2+|L(G)|-|L_2(G)|$. As such, $\beta(\mathcal{B}(G))=|G|^2+|L(G)|-(|G|^2+|L(G)|-|L_2(G)|)=|L_2(G)|$.
			\item Domatic partition of $V(\mathcal{B}(G))$ does not exist since   $L(G)$ is the only dominating set in $\mathcal{B}(G)$.    
		\end{enumerate}
	\end{remark}

	\section{Relations between \, $\mathcal{B}(G)$ \, and  probabilities associated to finite groups}
	In this section, we obtain relations between $\mathcal{B}(G)$ and certain probabilities associated to finite groups.     Let $G$ be a finite group and $H$ be any given subgroup of $G$. The probability that a randomly chosen pair of elements of $G$ generate $H$ is called the probability generating  a given subgroup. We write $\Pr_H(G)$ to denote this probability. Therefore, 
	\begin{equation}\label{SGP}
		{\Pr}_H(G)= \frac{|\{(a, b) \in G \times G : \langle a, b \rangle = H\}|}{|G|^2}.
	\end{equation}
	\begin{obj}
		\begin{enumerate}
			\item $\Pr_H(G)= 1$ if and only if $H=G=\{1\}$.
			\item $\Pr_H(G)= 1-\frac{1}{|G|^2}$ if and only if $H=G$ is a group of prime order.
			\item $\Pr_H(G)= \frac{1}{|G|^2}$ if and only if $H=\{1\}$.
		\end{enumerate}
	\end{obj}
	\noindent Note that $\Pr_G(G) := \varphi_2(G)$ is the probability that a randomly chosen pair of elements of $G$ generate $G$. Dixon \cite{Di69} obtained a lower bound for $\Pr_{A_5}(A_5)$ for the first time. Results on $\Pr_G(G)$ for symmetric group and finite simple groups can be found in \cite{Ba89, LS95, LS96}. The study of the generalized version of $\varphi_2(G)$, viz.
	\[
	\varphi_n(G) = \frac{|\{(x_1,  \dots, x_n) \in G \times  \cdots \times G : \langle x_1, \dots, x_n\rangle = G\}|}{|G|^n} 
	\]
	goes back to Hall \cite{Hall36}. Results on $\varphi_n(G)$ can be found in \cite{Pak99}.

	The probability that a randomly chosen pair of elements of $G$ commute is called the commuting probability of $G$. It is also known as commutativity degree of $G$. We write $\Pr(G)$ to denote this probability. Therefore, 
	\begin{align*}
		\Pr(G) &= \frac{|\{(a, b) \in G \times G : ab = ba \}|}{|G|^2} \\
		&= \frac{|\{(a, b) \in G \times G : \langle a, b \rangle \text{ is abelian} \}|}{|G|^2}.
	\end{align*}
	The origin of $\Pr(G)$ lies in a paper of Erd$\ddot{\rm{o}}$s and Tur$\acute{\rm a}$n \cite{ET69}. Results on $\Pr(G)$ can be found in the survey \cite{DNP-13}. A relation between the number of edges in  commuting/non-commuting graph and $\Pr(G)$ of $G$ was observed in \cite{AKM06, TE-13}.   Notions similar to $\Pr(G)$, viz. cyclicity degree (denoted by $\Pr_{\cyc}(G)$ and introduced in \cite{PSSW93}), nilpotency degree (denoted by $\Pr_{\nil}(G)$ and introduced in \cite{DGMW92}) and solvability degree (denoted by ${\Pr}_{\sol}(G)$ and introduced in \cite{FGSV2000})  are defined as follows:
	\[
	{\Pr}_{\cyc}(G)= \frac{|\{(a, b) \in G \times G : \langle a, b \rangle \text{ is cyclic} \}|}{|G|^2},
	\]
	\[
	{\Pr}_{\nil}(G)= \frac{|\{(a, b) \in G \times G : \langle a, b \rangle \text{ is nilpotent} \}|}{|G|^2}
	\]
	and 
	\[
	{\Pr}_{\sol}(G)= \frac{|\{(a, b) \in G \times G : \langle a, b \rangle \text{ is solvable} \}|}{|G|^2}.
	\]
	Relation between the number of edges in  solvable/non-solvable graph and ${\Pr}_{\sol}(G)$ of $G$ was observed in \cite{BNN2020}. Relations similar to \cite[Lemma 3.27]{AKM06} and \cite[Theorem 4.5]{BNN2020} can also be determined  for cyclic graph and nilpotent graph.
	
	In this section, we obtain certain relations among $\mathcal{B}(G)$, $\Pr(G)$, $\Pr_{\cyc}(G)$, $\Pr_{\nil}(G)$ and $\Pr_{\sol}(G)$. The following lemma is useful in this regard.
	\begin{lemma}\label{deg(H in L(G))}
		Let $G$ be a finite group and $H$ be a subgroup of $G$. Then 
		\[
		{\Pr}_H(G) = \frac{\deg_{\mathcal{B}(G)}(H)}{|G|^2}.
		\]
	\end{lemma}
	\begin{proof}
		We have  $\deg_{\mathcal{B}(G)}(H)=|\Nbd_{\mathcal{B}(G)}(H)|$, where 
		\begin{align*}
			\Nbd_{\mathcal{B}(G)}(H) & =\{(a, b) \in G \times G: (a, b) \text{ is adjacent to H}\}\\
			& =\{(a, b) \in G \times G: \langle a, b\rangle = H\}.
		\end{align*}
		Hence, the result follows from \eqref{SGP}.
	\end{proof}
	
	\begin{theorem}
		Let $G$ be a finite group and $e(\mathcal{B}(G))$ be the set of edges of the graph $\mathcal{B}(G)$. Then
		\[
		\sum_{(a, b)\in G\times G} \deg_{\mathcal{B}(G)}((a, b))= \sum_{H \in L(G)} \deg_{\mathcal{B}(G)}(H)= |G|^2=|e(\mathcal{B}(G))|
		\]
	\end{theorem}
	\begin{proof}
		For any bipartite graph $\mathcal{G}$ with partitions $A$ and $B$ of $V(\mathcal{G})$, we have
		\begin{equation}\label{deg_sum=num_of_edges}
			\sum_{x\in A} \deg_{\mathcal{G}}(x)= \sum_{y \in B} \deg_{\mathcal{G}}(y)=|e(\mathcal{G})|.
		\end{equation}
		Therefore, for the graph $\mathcal{B}(G)$ we have
		\[
		\sum_{(a, b)\in G\times G} \deg_{\mathcal{B}(G)}((a, b))= \sum_{H \in L(G)} \deg_{\mathcal{B}(G)}(H)= |e(\mathcal{B}(G))|.
		\]
		Since $\sum_{H \in L(G)}{\Pr}_H(G) = 1$, by Lemma \ref{deg(H in L(G))}, we have 
		\[
		\sum_{H \in L(G)} \deg_{\mathcal{B}(G)}(H)= \sum_{H \in L(G)}|G|^2 {\Pr}_H(G) = |G|^2.
		\]
		Hence the result follows.
	\end{proof}

	\begin{theorem}\label{relation between B(G) and varphi_2(G)}
		Let $G$ be a finite group and $H \in L(G)$. 
		Then 
		$
		\varphi_2(H) = \dfrac{\deg_{\mathcal{B}(G)}(H)}{|H|^2}.
		$
	\end{theorem}
	\begin{proof}
		For $a \neq b$,  the pairs $(a,b)$ and $(b,a)$ that generate $H$, contribute   one edge in $\Gamma_{\gen}(H)$ and two edges in $\mathcal{B}(G)$. It follows that
		\begin{align*}
			2|e(\Gamma_{\gen}(H))|
			&=\left|\{(a,b) \in H \times H: \langle a,b \rangle =H, a \neq b\}\right| \\
			&=\left|\{(a,b) \in H \times H: \langle a,b \rangle =H\}\right| - \left|\{(a,a) \in H \times H: \langle a \rangle=H\}\right| \\
			&= |H|^2 \varphi_2(H)-\phi(|H|),
		\end{align*}
		noting that $\varphi_2(H)=\frac{|\{(a,b) \in H \times H: \langle a, b \rangle=H\}|}{|H|^2}$ and $\phi(|H|) = 0$ or the number of generators of $\mathbb{Z}_{|H|}$ according as $H$ is non-cyclic or cyclic. Thus,
		\[
		|H|^2 \varphi_2(H) = 2|e(\Gamma_{\gen}(H))| + \phi(|H|).
		\]
		Hence, the result follows from Theorem \ref{relatn B(G) and generating graph}.
	\end{proof}

	\begin{theorem}
		Let $G$ be a finite group and $L_A(G)=\{H \in L(G): H \text{ is abelian}\}$. If $S=(G \times G) \sqcup L_{A}(G)$, then 
		\[
		\Pr(G) = \frac{\sum_{H \in L_A(G)} \deg_{\mathcal{B}(G)[S]}(H)}{|G|^2}=\frac{|e(\mathcal{B}(G)[S])|}{|G|^2}.
		\]
	\end{theorem}
	\begin{proof}
		Since $L_A(G)$ is the set of all abelian subgroups of $G$,  from the definitions of $\Pr_H(G)$ and $\Pr(G)$, we have
		\[
		\Pr(G)=\sum_{H \in L_A(G)}{\Pr}_H(G).
		\]
		Also, $\deg_{\mathcal{B}(G)[S]}(H) = \deg_{\mathcal{B}(G)}(H)$ for all $H \in L_A(G)$. Therefore,  using Lemma \ref{deg(H in L(G))}, we get
		\begin{align*}
			\sum_{H \in L_A(G)} \deg_{\mathcal{B}(G)[S]}(H)&=\sum_{H \in L_A(G)}|G|^2 {\Pr}_H(G) \\
			&=|G|^2\sum_{H \in L_A(G)}{\Pr}_H(G) \\
			&=|G|^2 \Pr(G).
		\end{align*}
		By \eqref{deg_sum=num_of_edges} we have 
		\[
		\sum_{H \in L_A(G)} \deg_{\mathcal{B}(G)[S]}(H)=|e(\mathcal{B}(G)[S])|.
		\]
		Hence, the result follows.
	\end{proof}
	\begin{theorem}
		Let $G$ be a finite group and $L_C(G)=\{H \in L(G): H \text{ is cyclic}\}$. If $S=(G \times G) \sqcup L_{C}(G)$, then 
		\[
		{\Pr}_{\cyc}(G)= \frac{\sum_{H \in L_C(G)} \deg_{\mathcal{B}(G)[S]}(H)}{|G|^2}= \frac{|e(\mathcal{B}(G)[S])|}{|G|^2}.
		\]
	\end{theorem}
	\begin{proof}
		Since $L_C(G)$ is the set of all cyclic subgroups of $G$,  from the definitions of $\Pr_H(G)$ and $\Pr_{\cyc}(G)$, we have
		\[
		{\Pr}_{\cyc}(G)=\sum_{H \in L_C(G)}{\Pr}_H(G).
		\]
		Also, $\deg_{\mathcal{B}(G)[S]}(H) = \deg_{\mathcal{B}(G)}(H)$ for all $H \in L_C(G)$. Therefore, using Lemma \ref{deg(H in L(G))}, we get
		\begin{align*}
			\sum_{H \in L_C(G)} \deg_{\mathcal{B}(G)[S]}(H)&=\sum_{H \in L_C(G)}|G|^2 {\Pr}_H(G) \\
			&=|G|^2\sum_{H \in L_C(G)}{\Pr}_H(G) \\
			&=|G|^2 {\Pr}_{\cyc}(G).
		\end{align*}
		By  \eqref{deg_sum=num_of_edges} we have 
		\[
		\sum_{H \in L_C(G)} \deg_{\mathcal{B}(G)[S]}(H)=|e(\mathcal{B}(G)[S])|.
		\]
		Hence, the result follows.
	\end{proof}
	\begin{theorem}
		Let $G$ be a finite group and $L_N(G)=\{H \in L(G): H \text{ is nilpotent}\}$. If $S=(G \times G) \sqcup L_{N}(G)$, then
		\[
		{\Pr}_{\nil}(G)= \frac{\sum_{H \in L_N(G)} \deg_{\mathcal{B}(G)[S]}(H)}{|G|^2}= \frac{|e(\mathcal{B}(G)[S])|}{|G|^2}.
		\]
	\end{theorem}
	\begin{proof}
		Since $L_N(G)$ is the set of all nilpotent subgroups of $G$,  from the definitions of $\Pr_H(G)$ and $\Pr_{\nil}(G)$, we have
		\[
		{\Pr}_{\nil}(G)=\sum_{H \in L_N(G)}{\Pr}_H(G).
		\]
		Also, $\deg_{\mathcal{B}(G)[S]}(H) = \deg_{\mathcal{B}(G)}(H)$ for all $H \in L_N(G)$. Therefore, using Lemma \ref{deg(H in L(G))}, we get
		\begin{align*}
			\sum_{H \in L_N(G)} \deg_{\mathcal{B}(G)[S]}(H)&=\sum_{H \in L_N(G)}|G|^2 {\Pr}_H(G) \\
			&=|G|^2\sum_{H \in L_N(G)}{\Pr}_H(G) \\
			&=|G|^2 {\Pr}_{\nil}(G).
		\end{align*}
		By \eqref{deg_sum=num_of_edges} we have 
		\[
		\sum_{H \in L_N(G)} \deg_{\mathcal{B}(G)[S]}(H)=|e(\mathcal{B}(G)[S])|.
		\]
		Hence, the result follows.
	\end{proof}

	\begin{theorem}
		Let $G$ be a finite group  and $L_S(G)=\{H \in L(G): H \text{ is solvable}\}$. If $S=(G \times G) \sqcup L_{S}(G)$, then
		\[
		{\Pr}_{\sol}(G)=\frac{\sum_{H \in L_S(G)} \deg_{\mathcal{B}(G)[S]}(H)}{|G|^2} = \frac{|e(\mathcal{B}(G)[S])|}{|G|^2}.
		\]
	\end{theorem}
	\begin{proof}
		Since $L_S(G)$ is the set of all solvable subgroups of $G$,  from the definitions of $\Pr_H(G)$ and $\Pr_{\sol}(G)$, we have
		\[
		{\Pr}_{\sol}(G)=\sum_{H \in L_S(G)}{\Pr}_H(G).
		\]
		Also, $\deg_{\mathcal{B}(G)[S]}(H) = \deg_{\mathcal{B}(G)}(H)$ for all $H \in L_S(G)$. Therefore, using Lemma \ref{deg(H in L(G))}, we get
		\begin{align*}
			\sum_{H \in L_S(G)} \deg_{\mathcal{B}(G)[S]}(H)&=\sum_{H \in L_S(G)}|G|^2 {\Pr}_H(G) \\
			&=|G|^2\sum_{H \in L_S(G)}{\Pr}_H(G) \\
			&=|G|^2 {\Pr}_{\sol}(G).
		\end{align*}
		By \eqref{deg_sum=num_of_edges} we have 
		$
		\sum_{H \in L_S(G)} \deg_{\mathcal{B}(G)[S]}(H)=|e(\mathcal{B}(G)[S])|.
		$
		Hence, the result follows.
	\end{proof}

	Let $\mathcal{G}$ be the set of all graphs.  A topological index  is a function $T : \mathcal{G} \to \mathbb{R}$ such that $T(\Gamma_1) = T(\Gamma_2)$ whenever the graphs $\Gamma_1$ and   $\Gamma_2$ are isomorphic.  
	Some of the well-known degree-based topological indices are Zagreb indices, Randic Connectivity index, Atom-Bond Connectivity index, Geometric-Arithmetic index, Harmonic index, Sum-Connectivity index etc. A survey on degree-based topological indices can be found in \cite{MNJ-FA-2020}. 
	Let $\Gamma \in \mathcal{G}$. The first and second Zagreb indices of $\Gamma$, denoted by $M_{1}(\Gamma)$ and $M_{2}(\Gamma)$ respectively, are defined as 
	\[
	M_{1}(\Gamma) = \sum\limits_{v \in V(\Gamma)} \deg(v)^{2}  \text{ and }  M_{2}(\Gamma) = \sum\limits_{uv \in e(\Gamma)} \deg(u)\deg(v).
	\]
	The Randic Connectivity index of $\Gamma$, denoted by $R(\Gamma)$, is defined as
	\[
	R(\Gamma)=\sum_{uv \in e(\Gamma)}\left(\deg(u)\deg(v)\right)^{\frac{-1}{2}}.
	\]
	The Atom-Bond Connectivity index of $\Gamma$, denoted by $\ABC(\Gamma)$, is defined as
	\[
	\ABC(\Gamma)=\sum_{uv\in e(\Gamma)}\left(\frac{\deg(u)+\deg(v)-2}{\deg(u)\deg(v)}\right)^{\frac{1}{2}}.
	\]
	The Geometric-Arithmetic index of $\Gamma$, denoted by $\GA(\Gamma)$, is defined as
	\[
	\GA(\Gamma)=\sum_{uv\in e(\Gamma)} \frac{\sqrt{\deg(u)\deg(v)}}{\frac{1}{2}(\deg(u)+\deg(v))}.
	\]
	The Harmonic index of $\Gamma$, denoted by $H(\Gamma)$, is defined as
	\[
	H(\Gamma)=\sum_{uv\in e(\Gamma)}\frac{2}{\deg(u)+\deg(v)}.
	\]
	The Sum-Connectivity index of $\Gamma$, denoted by $\SCI(\Gamma)$, is defined as
	\[
	\SCI(\Gamma)=\sum_{uv \in e(\Gamma)} \left(\deg(u)+\deg(v)\right)^{\frac{-1}{2}}.
	\]
	
	In the following theorem we obtain the above mentioned topological indices of $\mathcal{B}(G)$  in terms of $\varphi_2(G)$ using Theorem \ref{relation between B(G) and varphi_2(G)}.
	
	\begin{theorem}
		For any finite group $G$ we have  the following:
		\begin{enumerate}
			\item $M_1(\mathcal{B}(G))=|G|^2+ \sum\limits_{H \in L(G)}|H|^4\left(\varphi_2(H)\right)^2$
			and   
			
			\qquad\qquad\qquad\qquad\qquad\qquad\qquad\qquad\,$M_2(\mathcal{B}(G))=\sum\limits_{H \in L(G)}|H|^4 \left(\varphi_2(H)\right)^2$.
			\item $R(\mathcal{B}(G))=\sum\limits_{H \in L(G)}|H|\left(\varphi_2(H)\right)^{\frac{1}{2}}$.
			\item $\ABC(\mathcal{B}(G))=\sum\limits_{H \in L(G)}|H|\left(\left(|H|\varphi_2(H)\right)^2-\varphi_2(H)\right)^{\frac{1}{2}}$.
			\item $\GA(\mathcal{B}(G))=\sum\limits_{H \in L(G)} \frac{2|H|^3\left( \varphi_2(H)\right)^{\frac{3}{2}}}{(1+|H|^2 \varphi_2(H))}$.
			\item $H(\mathcal{B}(G)=\sum\limits_{H \in L(G)}\frac{2|H|^2 \varphi_2(H)}{1+|H|^2 \varphi_2(H)}$.
			\item $\SCI(\mathcal{B}(G))=\sum\limits_{H \in L(G)} |H|^2 \varphi_2(H)\left(1+|H|^2 \varphi_2(H)\right)^{\frac{-1}{2}}$.
		\end{enumerate}
	\end{theorem}
	\begin{proof}
		For $(a, b) \in G \times G$ and $H \in L(G)$, by Observation \ref{vrtex_deg_of_X_in_B(G)}(a) and Theorem \ref{relation between B(G) and varphi_2(G)}, we have $\deg_{\mathcal{B}(G)}((a, b))=1$ and $\deg_{\mathcal{B}(G)}(H)=|H|^2 \varphi_2(H)$. 
		\begin{enumerate}
			\item We have
			\begin{align*}
				M_{1}(\mathcal{B}(G)) &= \sum_{v \in V(\mathcal{B}(G))} \deg(v)^{2} \\
				&=\sum_{(a, b) \in G \times G}\left(\deg_{\mathcal{B}(G)}((a, b))\right)^2+\sum_{H \in L(G)}\left(\deg_{\mathcal{B}(G)}(H)\right)^2\\
				&= \sum_{(a, b) \in G \times G} 1 + \sum_{H \in L(G)}\left(|H|^2\varphi_2(H)\right)^2 
				=|G|^2+ \sum_{H \in L(G)}|H|^4\left(\varphi_2(H)\right)^2.
			\end{align*}
			Also, 
			\begin{align*}
				M_2&(\mathcal{B}(G))= \sum_{uv \in e(\mathcal{B}(G))} \deg(u)\deg(v)\\
				&\quad=\sum_{(a, b)H \in e(\mathcal{B}(G))}\deg_{\mathcal{B}(G)}((a, b))\deg_{\mathcal{B}(G)}(H) 
				= \sum_{(a, b)H \in e(\mathcal{B}(G))} \deg_{\mathcal{B}(G)}(H). \end{align*}
			In the above sum, $\deg_{\mathcal{B}(G)}(H)$ appears $\deg_{\mathcal{B}(G)}(H)$ many times for each $H \in L(G)$.     
			Therefore,      
			\begin{align*}
				M_2(\mathcal{B}(G))&= \sum_{H \in L(G)} \left(\deg_{\mathcal{B}(G)}(H)\right)^2 \\
				&= \sum_{H \in L(G)}\left(|H|^2 \varphi_2(H)\right)^2 
				= \sum_{H \in L(G)}|H|^4 \left(\varphi_2(H)\right)^2.
			\end{align*}
			\item We have
			\begin{align*}
				R(\mathcal{B}(G))&=\sum_{uv \in e(\mathcal{B}(G))}\left(\deg(u)\deg(v)\right)^{\frac{-1}{2}} \\
				&= \sum_{(a, b)H \in e(\mathcal{B}(G))}\left(\deg_{\mathcal{B}(G)}((a, b))\deg_{\mathcal{B}(G)}(H)\right)^{\frac{-1}{2}} \\
				&= \sum_{(a, b)H \in e(\mathcal{B}(G))}\left(\deg_{\mathcal{B}(G)}(H)\right)^{\frac{-1}{2}}.
			\end{align*}         
			In the above sum, \quad $\left(\deg_{\mathcal{B}(G)}(H)\right)^{\frac{-1}{2}}$ appears $\deg_{\mathcal{B}(G)}(H)$ many times for each $H \in L(G)$. Therefore,
			\begin{align*}          
				R&(\mathcal{B}(G))= \sum_{H \in L(G)}\left(\deg_{\mathcal{B}(G)}(H)\right)^{\frac{-1}{2}} \deg_{\mathcal{B}(G)}(H) \\
				&= \sum_{H \in L(G)}\left(\deg_{\mathcal{B}(G)}(H)\right)^{\frac{1}{2}} = \sum_{H \in L(G)}\left(|H|^2 \varphi_2(H)\right)^{\frac{1}{2}}=\sum_{H \in L(G)}|H|\left(\varphi_2(H)\right)^{\frac{1}{2}}.
			\end{align*}
			\item We have
			\begin{align*}
				\ABC(\mathcal{B}(G))&=\sum_{(a, b)H \in e(\mathcal{B}(G))}\left(\frac{\deg_{\mathcal{B}(G)}((a, b))+\deg_{\mathcal{B}(G)}(H)-2}{\deg_{\mathcal{B}(G)}((a, b))\deg_{\mathcal{B}(G)}(H)}\right)^{\frac{1}{2}} \\
				&=\sum_{(a, b)H \in e(\mathcal{B}(G))}\left(\frac{1+\deg_{\mathcal{B}(G)}(H)-2}{\deg_{\mathcal{B}(G)}(H)}\right)^{\frac{1}{2}} \\
				&=\sum_{(a, b)H \in e(\mathcal{B}(G))}\left(\frac{\deg_{\mathcal{B}(G)}(H)-1}{\deg_{\mathcal{B}(G)}(H)}\right)^{\frac{1}{2}}.
			\end{align*}
			In the above sum, \quad $\left(\frac{\deg_{\mathcal{B}(G)}(H)-1}{\deg_{\mathcal{B}(G)}(H)}\right)^{\frac{1}{2}}$ appears $\deg_{\mathcal{B}(G)}(H)$ many times for each $H \in L(G)$. Therefore,           
			\begin{align*}
				\ABC(\mathcal{B}(G)) &=\sum_{H \in L(G)}\left(\frac{\deg_{\mathcal{B}(G)}(H)-1}{\deg_{\mathcal{B}(G)}(H)}\right)^{\frac{1}{2}} \deg_{\mathcal{B}(G)}(H) \\
				&= \sum_{H \in L(G)}\left(\left(\deg_{\mathcal{B}(G)}(H)\right)^2-\deg_{\mathcal{B}(G)}(H)\right)^{\frac{1}{2}} \\
				&=\sum_{H \in L(G)}\left(\left(|H|^2 \varphi_2(H)\right)^2-|H|^2 \varphi_2(H)\right)^{\frac{1}{2}}.
			\end{align*}
			Hence, the result follows.
			\item We have
			\begin{align*}
				\GA(\mathcal{B}(G))&=\sum_{(a, b)H \in e(\mathcal{B}(G))} \frac{\sqrt{\deg_{\mathcal{B}(G)}((a, b))\deg_{\mathcal{B}(G)}(H)}}{\frac{1}{2}(\deg_{\mathcal{B}(G)}((a, b))+\deg_{\mathcal{B}(G)}(H))} \\
				&=\sum_{(a, b)H \in e(\mathcal{B}(G))} \frac{2\sqrt{\deg_{\mathcal{B}(G)}(H)}}{(1+\deg_{\mathcal{B}(G)}(H))}
			\end{align*}             
			In the above sum, \quad $\frac{2\sqrt{\deg_{\mathcal{B}(G)}(H)}}{(1+\deg_{\mathcal{B}(G)}(H))}$ \quad appears $\deg_{\mathcal{B}(G)}(H)$ many times for each $H \in L(G)$. Therefore,               
			
			\begin{align*}
				\ABC(\mathcal{B}(G))             
				&=\sum_{H \in L(G)} \frac{2\sqrt{\deg_{\mathcal{B}(G)}(H)}}{(1+\deg_{\mathcal{B}(G)}(H))}\deg_{\mathcal{B}(G)}(H)\\
				&= \sum_{H \in L(G)} \frac{2\left(\deg_{\mathcal{B}(G)}(H)\right)^{\frac{3}{2}}}{(1+\deg_{\mathcal{B}(G)}(H))} = \sum_{H \in L(G)} \frac{2\left(|H|^2 \varphi_2(H)\right)^{\frac{3}{2}}}{(1+|H|^2 \varphi_2(H))}.
			\end{align*}
			Hence, the result follows.
			\item We have
			\begin{align*}
				H(\mathcal{B}(G))&=\sum_{(a, b)H \in e(\mathcal{B}(G))}\frac{2}{\deg_{\mathcal{B}(G)}((a, b))+\deg_{\mathcal{B}(G)}(H)} \\
				&= \sum_{(a, b)H \in e(\mathcal{B}(G))}\frac{2}{1+\deg_{\mathcal{B}(G)}(H)}
			\end{align*}
			In the above sum, $\frac{2}{1+\deg_{\mathcal{B}(G)}(H)}$ appears $\deg_{\mathcal{B}(G)}(H)$ many times for each $H \in L(G)$. Therefore,          
			\begin{align*}
				H(\mathcal{B}(G))
				&=\sum_{H \in L(G)}\frac{2}{1+\deg_{\mathcal{B}(G)}(H)}\deg_{\mathcal{B}(G)}(H) 
				=\sum_{H \in L(G)}\frac{2|H|^2 \varphi_2(H)}{1+|H|^2 \varphi_2(H)}.
			\end{align*}
			\item We have
			\begin{align*}
				\SCI(\mathcal{B}(G))&=\sum_{(a, b)H \in e(\mathcal{B}(G))} \left(\deg_{\mathcal{B}(G)}((a, b))+\deg_{\mathcal{B}(G)}(H)\right)^{\frac{-1}{2}} \\
				&= \sum_{(a, b)H \in e(\mathcal{B}(G))} \left(1+\deg_{\mathcal{B}(G)}(H)\right)^{\frac{-1}{2}}.
			\end{align*}      
			In the above sum, $\left(1+\deg_{\mathcal{B}(G)}(H)\right)^{\frac{-1}{2}}$ appears $\deg_{\mathcal{B}(G)}(H)$ many times for each $H \in L(G)$. Therefore,       
			\begin{align*}
				\SCI(\mathcal{B}(G))&=\sum_{H \in L(G)} \left(1+\deg_{\mathcal{B}(G)}(H)\right)^{\frac{-1}{2}}\deg_{\mathcal{B}(G)}(H) \\
				&=\sum_{H \in L(G)} |H|^2 \varphi_2(H)\left(1+|H|^2 \varphi_2(H)\right)^{\frac{-1}{2}}.
			\end{align*}
		\end{enumerate}
\vspace{-.8cm}
	\end{proof}

	We conclude this section with the following  table describing  $|L_C(G)|$, $|L_A(G)|$, $|L_N(G)|$, $|L_S(G)|$, $|e(\mathcal{B}(G))|$ and various probabilities defined on finite groups for certain small order groups.
	\begin{table}[h]
		\begin{center}
			{{	
					\begin{tabular}{|c|c|c|c|c|c|c|c|}
						\hline
						$G$& $S_3$ & $D_8$ & $Q_8$ & $D_{10}$ & $D_{12}$ & $A_4$ & $S_4$ \\
						\hline
						$|G|$ & 6 & 8 & 8 & 10 & 12 & 12 & 24 \\
						\hline
						$|L_C(G)|$ & 5 & 7 & 5 & 7 & 10 & 8 & 17 \\
						\hline
						$|L_A(G)|$ & 5 & 9 & 5 & 7 & 13 & 9 & 21 \\
						\hline
						$|L_N(G)|$ & 5 & 10 & 6 & 7 & 13 & 9 & 24 \\
						\hline
						$|L_S(G)|$ & 6 & 10 & 6 & 8 & 16 & 10 & 30  \\
						\hline
						$|e(\mathcal{B}(G))|$ & 36 & 64 & 64 & 100 & 144 & 144 & 576 \\
						\hline
						$\Pr_{\cyc}(G)$ & $\frac{1}{2}$ & $\frac{7}{16}$ & $\frac{5}{8}$ & $\frac{2}{5}$ &$ \frac{3}{8}$ & $\frac{7}{24}$ & $\frac{1}{6}$  \\
						\hline
						$\Pr(G)$ & $\frac{1}{2}$ & $\frac{5}{8}$ & $\frac{5}{8}$ & $\frac{2}{5}$ & $\frac{1}{2}$ & $\frac{1}{3}$ & $\frac{5}{24}$ \\
						\hline
						$\Pr_{\nil}(G)$ & $\frac{1}{2}$ & 1 & 1 & $\frac{2}{5}$ & $\frac{1}{2}$ & $\frac{1}{3}$ & $\frac{1}{3}$ \\
						\hline
						$\Pr_{\sol}(G)$ & 1 & 1 & 1 & 1 & 1 & 1 & 1 \\
						\hline
						$\varphi_2(G)$ &$\frac{1}{2}$ & $\frac{3}{8}$ & $\frac{3}{8}$ & $\frac{3}{5}$ & $\frac{3}{8}$ & $\frac{2}{3}$ & $\frac{3}{8}$ \\
						\hline
					\end{tabular}
			}}
			\caption{Various probabilities of small order groups}\label{Table 1}
		\end{center}
	\end{table}	
	
	\newpage
	
	\section{Realization of $\mathcal{B}(G)$}
	Graph realization is one of the major aspects in studying graphs defined on algebraic systems.
	In Table \ref{Table 1}, while computing $|e(\mathcal{B}(G))|$ for various groups we realized  the structures of $\mathcal{B}(G)$ for $G = S_3, D_8, Q_8, D_{10}, D_{12}, A_4$ and $S_4$. For instance, $V(\mathcal{B}(S_3)) = S_3 \times S_3 \sqcup \{H_0, H_1, \dots, H_4, S_3\}$ where $H_0=\{(1)\}$, $H_1=\{(1), (12)\}$, $H_2=\{(1), (13)\}$, $H_3=\{(1), (23)\}$ and $H_4=\{(1), (123), (132)\}$. We have $\Nbd_{\mathcal{B}(S_3)}(H_0)=\{((1),(1))\}$, $\Nbd_{\mathcal{B}(S_3)}(H_i)=H_i \times H_i \setminus \{((1),(1))\}$ for $1 \leq i \leq 4$ and $\Nbd_{\mathcal{B}(S_3)}(S_3)=S_3 \times S_3 \setminus \left(\sqcup_{i=0}^{4} \Nbd_{\mathcal{B}(S_3)}(H_i)\right)$. Since the vertices from $S_3 \times S_3$ have degree one,  we have the following structure of $\mathcal{B}(S_3)$.
	\begin{center}
		\begin{tikzpicture}
			\tikzstyle{vertex}=[circle,minimum size=0.1pt,fill=black!30,inner sep=1.5pt]
			\node[vertex](A) at (-9.7,0){};
			\node[vertex](B) at (-9.7,-1){$H_0$};

			\node[vertex](C) at (-7.9,-0.8){$H_1$};
			\node[vertex](D) at (-7.9,0){};
			\node[vertex](E) at (-8.5,-1.6){};
			\node[vertex](F) at (-7.3,-1.6){};

			\node[vertex](C1) at (-6.1,-0.8){$H_2$};
			\node[vertex](D1) at (-6.1,0){};
			\node[vertex](E1) at (-6.7,-1.6){};
			\node[vertex](F1) at (-5.5,-1.6){};
			
			\node[vertex](C2) at (-4.3,-0.8){$H_3$};
			\node[vertex](D2) at (-4.3,0){};
			\node[vertex](E2) at (-4.9,-1.6){};
			\node[vertex](F2) at (-3.7,-1.6){};

			\node[vertex](G) at (-2,-0.8){$H_4$};
			\node[vertex](H) at (-2,0){};
			\node[vertex](I) at (-2,-1.6){};
			\node[vertex](J) at (-2.8,-0.8){};
			\node[vertex](K) at (-1.2,-0.8){};
			\node[vertex](L) at (-2.7,-0.2){};
			\node[vertex](M) at (-1.3,-0.2){};
			\node[vertex](N) at (-2.7,-1.4){};
			\node[vertex](O) at (-1.3,-1.4){};

			\node[vertex](P) at (0.8,-0.8){$S_3$};
			\node[vertex](Q) at (0.8,0){};
			\node[vertex](R) at (0.3,0){};
			\node[vertex](S) at (-0.2,0){};
			\node[vertex](T) at (1.3,0){};
			\node[vertex](U) at (1.8,0){};
			\node[vertex](V) at (2.1,-0.3){};
			\node[vertex](W) at (2.1,-0.6){};
			\node[vertex](X) at (2.1,-0.9){};
			\node[vertex](Y) at (2.1,-1.2){};
			\node[vertex](Z) at (1.8,-1.6){};
			\node[vertex](a) at (1.3,-1.6){};
			\node[vertex](b) at (0.8,-1.6){};
			\node[vertex](c) at (0.3,-1.6){};
			\node[vertex](d) at (-0.2,-1.6){};
			\node[vertex](e) at (-0.5,-0.3){};
			\node[vertex](f) at (-0.5,-0.6){};
			\node[vertex](g) at (-0.5,-0.9){};
			\node[vertex](h) at (-0.5,-1.2){};
			\path 
			(A) edge (B)
			
			(C) edge (D)
			(C) edge (E)
			(C) edge (F)
			
			(C1) edge (D1)
			(C1) edge (E1)
			(C1) edge (F1)
			
			(C2) edge (D2)
			(C2) edge (E2)
			(C2) edge (F2)
			
			(G) edge (H)
			(G) edge (K)
			(G) edge (J)
			(G) edge (I)
			(G) edge (L)
			(G) edge (M)
			(G) edge (N)
			(G) edge (O)

			(P) edge (Q)
			(P) edge (R)
			(P) edge (S)
			(P) edge (T)
			(P) edge (U)
			(P) edge (V)
			(P) edge (W)
			(P) edge (X)
			(P) edge (Y)
			(P) edge (Z)
			(P) edge (a)
			(P) edge (b)
			(P) edge (c)
			(P) edge (d)
			(P) edge (e)
			(P) edge (f)
			(P) edge (g)
			(P) edge (h);
		\end{tikzpicture}
		
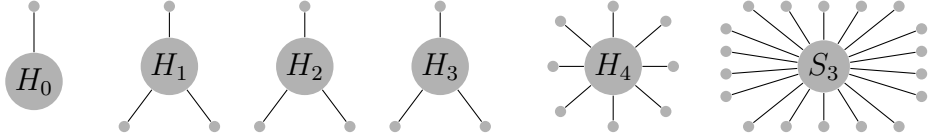
\captionof{figure}{Graph structure of $\mathcal{B}(S_3)$}
		\label{fig:fig1}
	\end{center}
	From Figure \ref{fig:fig1}, it is clear that $\mathcal{B}(S_3)=K_2 \sqcup 3K_{1, 3} \sqcup K_{1, 8} \sqcup K_{1, 18}$. We know that $D_6 \cong S_3$. Therefore, $\mathcal{B}(D_6) \cong \mathcal{B}(S_3) =K_2 \sqcup 3K_{1, 3} \sqcup K_{1, 8} \sqcup K_{1, 18}$.
	
	For the group $Q_8$, we have $V(\mathcal{B}(Q_8)) =Q_8\times Q_8 \sqcup \{H_0, H_1, \dots, H_4, Q_8\}$ where $H_0=\{1\}$, $H_1=\{1, a^2\}$, $H_2=\{1, a, a^2, a^3\}$, $H_3=\{1, a^2, b, a^2b\}$ and $H_4=\{1, a^2, ab, a^3b\}$. We have $\Nbd_{\mathcal{B}(Q_8)}(H_0)=\{(1,1)\}$, $\Nbd_{\mathcal{B}(Q_8)}(H_1)=H_1 \times H_1 \setminus \{(1,1)\}$,  $\Nbd_{\mathcal{B}(Q_8)}(H_i)=H_i \times H_i \setminus \{\sqcup_{j=0}^{1}\Nbd_{\mathcal{B}(Q_8)}(H_j)\}$ for $2 \leq i \leq 4$ and $\Nbd_{\mathcal{B}(Q_8)}(Q_8)=Q_8 \times Q_8 \setminus \left(\sqcup_{i=0}^{4} \Nbd_{\mathcal{B}(Q_8)}(H_i)\right)$. Since the vertices from $Q_8 \times Q_8$ have degree one, we have the following structure of $\mathcal{B}(Q_8)$.
	\begin{center}
		\begin{tikzpicture}
			\tikzstyle{vertex}=[circle,minimum size=0.1pt,fill=black!30,inner sep=1.5pt]
			\node[vertex](A) at (-8.6,0){};
			\node[vertex](B) at (-8.6,-1){$H_0$};

			\node[vertex](C) at (-6,-0.8){$H_1$};
			\node[vertex](D) at (-6,0){};
			\node[vertex](E) at (-6.8,-1.6){};
			\node[vertex](F) at (-5.2,-1.6){};

			\node[vertex](C2) at (-3,-0.8){$H_2$};
			\node[vertex](D2) at (-3,0){};
			\node[vertex](E2) at (-3,-1.6){};
			\node[vertex](F2) at (-3.9,-0.8){};
			\node[vertex](G2) at (-2.1,-0.8){};
			\node[vertex](H2) at (-3.5,-0.1){};
			\node[vertex](I2) at (-3.8,-0.4){};
			\node[vertex](J2) at (-2.5,-0.1){};
			\node[vertex](K2) at (-2.2,-0.4){};
			\node[vertex](L2) at (-3.8,-1.2){};
			\node[vertex](M2) at (-3.5,-1.5){};
			\node[vertex](N2) at (-2.2,-1.2){};
			\node[vertex](O2) at (-2.5,-1.5){};

			\node[vertex](G) at (0.5,-0.8){$H_3$};
			\node[vertex](H) at (0.5,0){};
			\node[vertex](I) at (0.5,-1.6){};
			\node[vertex](J) at (-0.4,-0.8){};
			\node[vertex](K) at (1.4,-0.8){};
			\node[vertex](L) at (0,-0.1){};
			\node[vertex](L1) at (-0.3,-0.4){};
			\node[vertex](M) at (1,-0.1){};
			\node[vertex](M1) at (1.3,-0.4){};
			\node[vertex](N) at (1.3,-1.2){};
			\node[vertex](N1) at (0,-1.5){};
			\node[vertex](O) at (-0.3,-1.2){};
			\node[vertex](O1) at (1,-1.5){};

			\node[vertex](C3) at (-6,-3.2){$H_4$};
			\node[vertex](D3) at (-6,-2.2){};
			\node[vertex](E3) at (-6,-4.2){};
			\node[vertex](F3) at (-7.1,-3.2){};
			\node[vertex](G3) at (-4.9,-3.2){};
			\node[vertex](H3) at (-6.6,-2.3){};
			\node[vertex](I3) at (-7,-2.7){};
			\node[vertex](J3) at (-5.4,-2.3){};
			\node[vertex](K3) at (-5,-2.7){};
			\node[vertex](L3) at (-7,-3.7){};
			\node[vertex](M3) at (-6.6,-4.1){};
			\node[vertex](N3) at (-5,-3.7){};
			\node[vertex](O3) at (-5.4,-4.1){};

			\node[vertex](P) at (-2,-3.2){$Q_8$};
			\node[vertex](Q) at (-2,-2.2){};
			\node[vertex](R) at (-2.5,-4.2){};
			\node[vertex](S) at (-3,-2.2){};
			\node[vertex](m) at (-3.35,-2.35){};
			\node[vertex](n) at (-0.65,-2.35){};
			\node[vertex](T) at (-1.5,-2.2){};
			\node[vertex](U) at (-1,-2.2){};
			\node[vertex](V) at (-0.5,-2.6){};
			\node[vertex](W) at (-0.5,-2.9){};
			\node[vertex](X) at (-0.5,-3.2){};
			\node[vertex](Y) at (-0.5,-3.5){};
			\node[vertex](k) at (-0.5,-3.8){};
			\node[vertex](l) at (-0.6,-4.1){};
			\node[vertex](Z) at (-1,-4.2){};
			\node[vertex](a) at (-1.5,-4.2){};
			\node[vertex](b) at (-2,-4.2){};
			\node[vertex](c) at (-2.5,-2.2){};
			\node[vertex](d) at (-3,-4.2){};
			\node[vertex](e) at (-3.5,-2.6){};
			\node[vertex](f) at (-3.5,-2.9){};
			\node[vertex](g) at (-3.5,-3.2){};
			\node[vertex](h) at (-3.5,-3.5){};
			\node[vertex](i) at (-3.5,-3.8){};
			\node[vertex](j) at (-3.35,-4.1){};
			\path 
			(A) edge (B)
			
			(C) edge (D)
			(C) edge (E)
			(C) edge (F)
			(C2) edge (D2)
			(C2) edge (E2)
			(C2) edge (F2)
			(C2) edge (G2)
			(C2) edge (H2)
			(C2) edge (I2)
			(C2) edge (J2)
			(C2) edge (K2)
			(C2) edge (L2)
			(C2) edge (M2)
			(C2) edge (N2)
			(C2) edge (O2)
			
			(G) edge (H)
			(G) edge (K)
			(G) edge (J)
			(G) edge (I)
			(G) edge (L)
			(G) edge (M)
			(G) edge (N)
			(G) edge (O)
			(G) edge (L1)
			(G) edge (M1)
			(G) edge (N1)
			(G) edge (O1)
			
			(C3) edge (D3)
			(C3) edge (E3)
			(C3) edge (F3)
			(C3) edge (G3)
			(C3) edge (H3)
			(C3) edge (I3)
			(C3) edge (J3)
			(C3) edge (K3)
			(C3) edge (L3)
			(C3) edge (M3)
			(C3) edge (N3)
			(C3) edge (O3)

			(P) edge (Q)
			(P) edge (R)
			(P) edge (S)
			(P) edge (T)
			(P) edge (U)
			(P) edge (V)
			(P) edge (W)
			(P) edge (X)
			(P) edge (Y)
			(P) edge (Z)
			(P) edge (a)
			(P) edge (b)
			(P) edge (c)
			(P) edge (d)
			(P) edge (e)
			(P) edge (f)
			(P) edge (g)
			(P) edge (h)
			(P) edge (i)
			(P) edge (j)
			(P) edge (k)
			(P) edge (l)
			(P) edge (m)
			(P) edge (n);
		\end{tikzpicture}
		
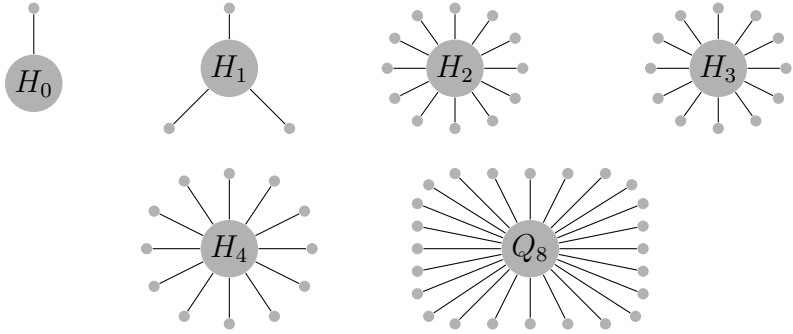
\captionof{figure}{Graph structure of $\mathcal{B}(Q_8)$}
		\label{fig:fig2}
	\end{center}
	From Figure \ref{fig:fig2}, it is clear that $\mathcal{B}(Q_8) = K_2 \sqcup K_{1, 3} \sqcup 3K_{1, 12} \sqcup K_{1, 24}$.
	
	For the group $D_8$, we have $V(\mathcal{B}(D_8))=D_8 \times D_8 \sqcup \{H_0, H_1, \dots, H_8, D_8\}$ where $H_0=\{1\}$, $H_1=\{1, a^2\}$, $H_2=\{1, b\}$, $H_3=\{1, ab\}$, $H_4=\{1, a^2b\}$, $H_5=\{1, a^3b\}$, $H_6=\{1, a^2, b, a^2b\}$, $H_7=\{1, a^2, ab, a^3b\}$ and $H_8=\{1, a, a^2, a^3\}$. We have $\Nbd_{\mathcal{B}(D_8)}(H_0)=\{(1,1)\}$, $\Nbd_{\mathcal{B}(D_8)}(H_i)=H_i \times H_i \setminus \{(1,1)\}$ for $1 \leq i \leq 5$, $\Nbd_{\mathcal{B}(D_8)}(H_6)= H_6 \times H_6 \setminus \left(\sqcup_{j=0}^{2}\Nbd_{\mathcal{B}(D_8)}(H_j) \sqcup \Nbd_{\mathcal{B}(D_8)}(H_4)\right)$, $\Nbd_{\mathcal{B}(D_8)}(H_7)= H_7 \times H_7 \setminus \left(\sqcup_{j=0}^{1}\Nbd_{\mathcal{B}(D_8)}(H_j)\right.$ $\left. \sqcup \Nbd_{\mathcal{B}(D_8)}(H_3) \sqcup \Nbd_{\mathcal{B}(D_8)}(H_5)\right)$, $\Nbd_{\mathcal{B}(D_8)}(H_8)=H_8 \times H_8 \setminus \left(\sqcup_{j=0}^{1}\Nbd_{\mathcal{B}(D_8)}(H_j)\right)$ and $\Nbd_{\mathcal{B}(D_8)}(D_8)=D_8 \times D_8$ $ \setminus \left(\sqcup_{j=0}^{8}\Nbd_{\mathcal{B}(D_8)}(H_j)\right)$. Since the vertices from $D_8 \times D_8$ have degree one, we have the following structure of $\mathcal{B}(D_8)$.
	\begin{center}
		\begin{tikzpicture}
			\tikzstyle{vertex}=[circle,minimum size=0.1pt,fill=black!30,inner sep=1.5pt]
			\node[vertex](A) at (-10.7,0){};
			\node[vertex](B) at (-10.7,-1){$H_0$};
			
			\node[vertex](C) at (-8.7,-0.8){$H_1$};
			\node[vertex](D) at (-8.7,0){};
			\node[vertex](E) at (-9.3,-1.6){};
			\node[vertex](F) at (-8.1,-1.6){};
			
			\node[vertex](C4) at (-6.5,-0.8){$H_2$};
			\node[vertex](D4) at (-6.5,0){};
			\node[vertex](E4) at (-7.1,-1.6){};
			\node[vertex](F4) at (-5.9,-1.6){};
			
			\node[vertex](C5) at (-4.2,-0.8){$H_3$};
			\node[vertex](D5) at (-4.2,0){};
			\node[vertex](E5) at (-4.8,-1.6){};
			\node[vertex](F5) at (-3.6,-1.6){};

			\node[vertex](C2) at (-2,-0.8){$H_4$};
			\node[vertex](D2) at (-2,0){};
			\node[vertex](E2) at (-2.6,-1.6){};
			\node[vertex](F2) at (-1.4,-1.6){};

			\node[vertex](G) at (0.2,-0.8){$H_5$};
			\node[vertex](H) at (0.2,0){};
			\node[vertex](I) at (-0.4,-1.6){};
			\node[vertex](J) at (0.8,-1.6){};

			\node[vertex](a1) at (-9.8,-3.2){$H_6$};
			\node[vertex](b1) at (-9.8,-2.2){};
			\node[vertex](c1) at (-9.8,-4.2){};
			\node[vertex](d1) at (-10.9,-2.7){};
			\node[vertex](e1) at (-10.9,-3.7){};
			\node[vertex](f1) at (-8.7,-2.7){};
			\node[vertex](g1) at (-8.7,-3.7){};

			\node[vertex](a2) at (-6.8,-3.2){$H_7$};
			\node[vertex](b2) at (-6.8,-2.2){};
			\node[vertex](c2) at (-6.8,-4.2){};
			\node[vertex](d2) at (-7.9,-2.7){};
			\node[vertex](e2) at (-7.9,-3.7){};
			\node[vertex](f2) at (-5.7,-2.7){};
			\node[vertex](g2) at (-5.7,-3.7){};

			\node[vertex](C3) at (-3.8,-3.2){$H_8$};
			\node[vertex](D3) at (-3.8,-2.2){};
			\node[vertex](E3) at (-3.8,-4.2){};
			\node[vertex](F3) at (-4.9,-3.2){};
			\node[vertex](G3) at (-2.7,-3.2){};
			\node[vertex](H3) at (-4.4,-2.3){};
			\node[vertex](I3) at (-4.8,-2.7){};
			\node[vertex](J3) at (-3.2,-2.3){};
			\node[vertex](K3) at (-2.8,-2.7){};
			\node[vertex](L3) at (-4.8,-3.7){};
			\node[vertex](M3) at (-4.4,-4.1){};
			\node[vertex](N3) at (-2.8,-3.7){};
			\node[vertex](O3) at (-3.2,-4.1){};

			\node[vertex](P) at (-0.6,-3.2){$D_8$};
			\node[vertex](Q) at (-0.6,-2.2){};
			\node[vertex](R) at (-1.1,-4.2){};
			\node[vertex](S) at (-1.6,-2.2){};
			\node[vertex](m) at (-1.95,-2.35){};
			\node[vertex](n) at (0.75,-2.35){};
			\node[vertex](T) at (-0.1,-2.2){};
			\node[vertex](U) at (0.4,-2.2){};
			\node[vertex](V) at (0.9,-2.6){};
			\node[vertex](W) at (0.9,-2.9){};
			\node[vertex](X) at (0.9,-3.2){};
			\node[vertex](Y) at (0.9,-3.5){};
			\node[vertex](k) at (0.9,-3.8){};
			\node[vertex](l) at (0.9,-4.1){};
			\node[vertex](Z) at (0.4,-4.2){};
			\node[vertex](a) at (-0.1,-4.2){};
			\node[vertex](b) at (-0.6,-4.2){};
			\node[vertex](c) at (-1.1,-2.2){};
			\node[vertex](d) at (-1.6,-4.2){};
			\node[vertex](e) at (-2.1,-2.6){};
			\node[vertex](f) at (-2.1,-2.9){};
			\node[vertex](g) at (-2.1,-3.2){};
			\node[vertex](h) at (-2.1,-3.5){};
			\node[vertex](i) at (-2.1,-3.8){};
			\node[vertex](j) at (-1.95,-4.1){};
			\path 
			(A) edge (B)
			
			(C) edge (D)
			(C) edge (E)
			(C) edge (F)
			
			(C2) edge (D2)
			(C2) edge (E2)
			(C2) edge (F2)
			
			(C4) edge (D4)
			(C4) edge (E4)
			(C4) edge (F4)
			
			(C5) edge (D5)
			(C5) edge (E5)
			(C5) edge (F5)

			(G) edge (H)
			(G) edge (I)
			(G) edge (J)
			
			(a1) edge (b1)
			(a1) edge (c1)
			(a1) edge (d1)
			(a1) edge (e1)
			(a1) edge (f1)
			(a1) edge (g1)

			(a2) edge (b2)
			(a2) edge (c2)
			(a2) edge (d2)
			(a2) edge (e2)
			(a2) edge (f2)
			(a2) edge (g2)

			(C3) edge (D3)
			(C3) edge (E3)
			(C3) edge (F3)
			(C3) edge (G3)
			(C3) edge (H3)
			(C3) edge (I3)
			(C3) edge (J3)
			(C3) edge (K3)
			(C3) edge (L3)
			(C3) edge (M3)
			(C3) edge (N3)
			(C3) edge (O3)

			(P) edge (Q)
			(P) edge (R)
			(P) edge (S)
			(P) edge (T)
			(P) edge (U)
			(P) edge (V)
			(P) edge (W)
			(P) edge (X)
			(P) edge (Y)
			(P) edge (Z)
			(P) edge (a)
			(P) edge (b)
			(P) edge (c)
			(P) edge (d)
			(P) edge (e)
			(P) edge (f)
			(P) edge (g)
			(P) edge (h)
			(P) edge (i)
			(P) edge (j)
			(P) edge (k)
			(P) edge (l)
			(P) edge (m)
			(P) edge (n);
		\end{tikzpicture}
		\captionof{figure}{Graph structure of $\mathcal{B}(D_8)$}
		\label{fig:fig3}
	\end{center}
	From Figure \ref{fig:fig3}, it is clear that $\mathcal{B}(D_8) = K_2 \sqcup 5K_{1, 3} \sqcup 2K_{1, 6} \sqcup K_{1, 12} \sqcup K_{1, 24}$. Thus, it follows that $\mathcal{B}(Q_8)$ and  $\mathcal{B}(D_8)$ are not isomorphic.  It is worth mentioning that the commuting and non-commuting graphs of the groups $Q_8$ and $D_8$ are isomorphic. 
	\begin{center}
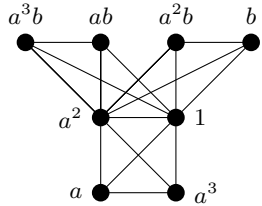

		\begin{tikzpicture}
			\Vertex[x=0, y=-3, size=0.05, color=black, label=$a$, position=left]{A}
			\Vertex[x=1, y=-3, size=0.05, color=black, label=$a^3$, position=right]{B}
			\Vertex[x=0, y=-2, size=0.05, color=black, label=$a^2$, position=left]{C}
			\Vertex[x=0, y=-1, size=0.05, color=black, label=$ab$, position=above]{D}
			\Vertex[x=-1, y=-1, size=0.05, color=black, label=$a^3b$, position=above]{E}
			\Vertex[x=1, y=-2, size=0.05, color=black, label=1, position=right]{F}
			\Vertex[x=1, y=-1, size=0.05, color=black, label=$a^2b$, position=above]{G}
			\Vertex[x=2, y=-1, size=0.05, color=black, label=$b$, position=above]{H}
			\path (A) edge (B)
			(A) edge (C)
			(A) edge (F)
			(C) edge (B)
			(F) edge (B)
			(C) edge (F)
			(C) edge (E)
			(C) edge (D)
			(C) edge (G)
			(C) edge (H)
			(F) edge (E)
			(F) edge (D)
			(F) edge (G)
			(F) edge (H)
			(C) edge (E)
			(C) edge (D)
			(C) edge (G)
			(D) edge (E)
			(C) edge (E)
			(C) edge (D)
			(C) edge (G)
			(G) edge (H);
		\end{tikzpicture}
		\captionof{figure}{Commuting graphs of $D_8$ and $Q_8$ }
		\label{fig:figA}
		
	\end{center}

	For the group $D_{10}$, we have $V(\mathcal{B}(D_{10}))= D_{10} \times D_{10} \sqcup \{H_0, H_1, \ldots, H_6, D_{10}\}$ where $H_0=\{1\}$, $H_1=\langle b \rangle$, $H_2=\langle ab \rangle$, $H_3=\langle a^2b \rangle$, $H_4=\langle a^3b \rangle$, $H_5=\langle a^4b \rangle$ and $H_6=\langle a \rangle$. We have $\Nbd_{\mathcal{B}(D_{10})}(H_0)=\{(1,1)\}$, $\Nbd_{\mathcal{B}(D_{10})}(H_i)=H_i \times H_i \setminus \{(1,1)\}$ for $1 \leq i \leq 6$ and $\Nbd_{\mathcal{B}(D_{10})}(D_{10})=D_{10} \times D_{10} \setminus \left(\sqcup_{i=0}^{6}\Nbd_{\mathcal{B}(D_{10})}(H_i)\right)$. Since the vertices from $D_{10} \times D_{10}$ have degree one, it follows that $\mathcal{B}(D_{10}) = K_2 \sqcup 5K_{1, 3} \sqcup K_{1, 24} \sqcup K_{1, 60}$.
	
	For the group $D_{12}$, we have $V(\mathcal{B}(D_{12}))= D_{12} \times D_{12} \sqcup \{H_0, H_1, \ldots, H_{14}, D_{12}\}$ where $H_0=\{1\}$, $H_1=\langle b \rangle$, $H_2=\langle ab \rangle$, $H_3=\langle a^2b \rangle$, $H_4=\langle a^3b \rangle$, $H_5=\langle a^4b \rangle$, $H_6=\langle a^5b \rangle$, $H_7=\langle a^3 \rangle$, $H_8= \langle a^2 \rangle$, $H_9=\langle a \rangle$, $H_{10}=\{1, a^3, b, a^3b\}$, $H_{11}=\{1, a^3, ab, a^4b\}$, $H_{12}=\{1, a^3, a^2b, a^5b\}$, $H_{13}=\{1, a^2, a^4, b, a^2b, a^4b\}$ and $H_{14}=\{1, a^2, a^4, ab, a^3b, a^5b\}$. We have $\Nbd_{\mathcal{B}(D_{12})}(H_0)=\{(1,1)\}$, $\Nbd_{\mathcal{B}(D_{12})}(H_i)=H_i \times H_i \setminus \{(1,1)\}$ for $1 \leq i \leq 8$, $\Nbd_{\mathcal{B}(D_{12})}(H_9)=H_9 \times H_9 \setminus (\sqcup_{j=7}^{8}\Nbd_{\mathcal{B}(D_{12})}(H_j)$ $ \sqcup \{(1,1)\})$ and $\Nbd_{\mathcal{B}(D_{12})}(H_{10})=H_{10} \times H_{10} \setminus (\Nbd_{\mathcal{B}(D_{12})}(H_1) \sqcup \Nbd_{\mathcal{B}(D_{12})}(H_4) \sqcup \Nbd_{\mathcal{B}(D_{12})}(H_7) \sqcup \{(1,1)\})$. Now, since $H_{10} \cong H_{11} \cong H_{12}$ and  $\mathcal{B}(D_{12})[\{H_{10}\} \sqcup \Nbd_{\mathcal{B}(D_{12})}(H_{10})] \cong K_{1, 6}$ so $\mathcal{B}(D_{12})[\{H_{i}\} \sqcup \Nbd_{\mathcal{B}(D_{12})}(H_i)] \cong K_{1, 6}$ for $i=11$ and $12$.
	Also, since $H_{13} \cong H_{14} \cong S_3$ so $\mathcal{B}(D_{12})[\{H_i\} \sqcup \Nbd_{\mathcal{B}(D_{12})}(H_i)] \cong K_{1, 18}$ for $i=13$ and $14$.
	Now, $\Nbd_{\mathcal{B}(D_{12})}(D_{12})=D_{12} \times D_{12} \setminus \left(\sqcup_{i=0}^{14}\Nbd_{\mathcal{B}(D_{12})}(H_i)\right)$. Since the vertices from $D_{12} \times D_{12}$ have degree one, it follows that $\mathcal{B}(D_{12}) = K_2 \sqcup 7K_{1, 3} \sqcup K_{1, 8} \sqcup K_{1, 24} \sqcup 3K_{1, 6} \sqcup 2K_{1, 18} \sqcup K_{1, 54}$.
	
	For the group $A_4$, we have $V(\mathcal{B}(A_4))=A_4 \times A_4 \sqcup \{H_0, H_1, \dots, H_8, A_4\}$ where $H_0=\{(1)\}$, $H_1=\{(1), (12)(34)\}$, $H_2=\{(1), (13)(24)\}$, $H_3=\{(1), (14)(23)\}$, $H_4=\{(1), (123)$ $, (132)\}$, $H_5=\{(1), (134), (143)\}$, $H_6=\{(1), (234), (243)\}$, $H_7=\{(1), (124), (142)\}$ and $H_8=\{(1), (12)(34), (13)(24), (14)(23)\}$. We have $\Nbd_{\mathcal{B}(A_4)}(H_0)=\{((1), (1))\}$, $\Nbd_{\mathcal{B}(A_4)}(H_i)=H_i \times H_i \setminus \{((1), (1))\}$ for $1 \leq i \leq 7$, $\Nbd_{\mathcal{B}(A_4)}(H_8)=H_8 \times H_8 \setminus \left(\sqcup_{j=0}^{3}\Nbd_{\mathcal{B}(A_4)}(H_j)\right)$ and $\Nbd_{\mathcal{B}(A_4)}(A_4)=A_4 \times A_4$ $ \setminus \left(\sqcup_{j=0}^{8}\Nbd_{\mathcal{B}(A_4)}(H_j)\right)$. Since the vertices from $A_4 \times A_4$ have degree one, it follows that $\mathcal{B}(A_4) = K_2 \sqcup 3K_{1, 3} \sqcup 4K_{1, 8} \sqcup K_{1, 6} \sqcup K_{1, 96}$.

	For the group $S_4$, we have $V(\mathcal{B}(S_4))= S_4 \times S_4 \sqcup \{H_0, H_1, \ldots, H_{28}, S_4\}$, where 
	\begin{center}
		$H_0=\{(1)\}$,
	\end{center}
	
	\begin{center}
		$\begin{array}{llll}
			H_1	=\langle (12) \rangle,  &\quad H_2=\langle (13) \rangle, & \quad H_3=\langle (14) \rangle, & \quad H_4=\langle (23) \rangle, \\
			H_5=\langle (24) \rangle, & \quad H_6 =\langle (34) \rangle,  &\quad H_7=\langle (12)(34) \rangle,  &\quad  H_8=\langle (13)(24) \rangle,\\
			H_9=\langle (14)(23) \rangle,	&\quad H_{10}=\langle (123) \rangle,  &\quad H_{11}=\langle (134) \rangle,  &\quad H_{12}=\langle (124) \rangle,  \\
			H_{13}=\langle (234) \rangle,	&\quad H_{14}=\langle (1234) \rangle,  &\quad H_{15}=\langle (1324) \rangle,  &\quad H_{16}=\langle (1243) \rangle, 
		\end{array}$
	\end{center}
	\begin{center}
		$\begin{array}{ll}
			H_{17}=\langle (12), (34) \rangle,	& \quad H_{18}=\langle (13), (24) \rangle,\\
			H_{19}=\langle (14), (23) \rangle,	& \quad H_{20}=\langle (12)(34), (13)(24) \rangle,
		\end{array}$
	\end{center}

	\noindent $H_{21}=\{(1), (12), (13), (23), (123), (132)\}$, $H_{22}=\{(1), (13), (14), (34), (134), (143)\}$, 
	
	\noindent $H_{23}=\{(1), (12), (14), (24), (124), (142)\}$, $H_{24}=\{(1), (23), (24), (34), (234), (242)\}$, 
	
	\noindent $H_{25}=\{(1), (12), (34), (12)(34), (13)(24), (14)(23), (1423), (1324)\}$, 
	
	\noindent $H_{26}=\{(1), (13), (24), (12)(34), (13)(24), (14)(23), (1234), (1432)\}$, 
	
	\noindent $H_{27}=\{(1), (14), (23), (12)(34), (13)(24), (14)(23), (1342), (1243)\}$ and 
	
	\noindent $H_{28}=$ 
	$\{(1), (12)(34), (13)(24), (14)(23)$, 
	
	\qquad\qquad\qquad\qquad\qquad\qquad$(123), (132), (134), (143), (234), (243), (124), (142)\}$.
	
	\noindent Note that the vertices from $S_4 \times S_4$ have degree one in $\mathcal{B}(S_4)$. We have $\Nbd_{\mathcal{B}(S_4)}(H_0)=\{((1), (1))\}$ and so $\mathcal{B}(S_4)[\{H_0\} \sqcup \Nbd_{\mathcal{B}(S_4)}(H_0)] = K_2$. For $1 \leq i \leq 9$, we have $H_i \cong \mathbb{Z}_2$ and 
	so $\mathcal{B}(S_4)[\{H_i\} \sqcup \Nbd_{\mathcal{B}(S_4)}(H_i)] = \mathcal{B}(\mathbb{Z}_2)[\{\mathbb{Z}_2\} \sqcup \Nbd_{\mathcal{B}(\mathbb{Z}_2)}(\mathbb{Z}_2)] = K_{1, 3}$.
	For  $10 \leq i \leq 13$, we have $H_i \cong \mathbb{Z}_3$ and so 
	$\mathcal{B}(S_4)[\{H_i\} \sqcup \Nbd_{\mathcal{B}(S_4)}(H_i)] = \mathcal{B}(\mathbb{Z}_3)[\{\mathbb{Z}_3\} \sqcup \Nbd_{\mathcal{B}(\mathbb{Z}_3)}(\mathbb{Z}_3)] = K_{1, 8}$.
	For $14 \leq i \leq 16$, we have $H_i \cong \mathbb{Z}_4$ and so $\mathcal{B}(S_4)[\{H_i\} \sqcup \Nbd_{\mathcal{B}(S_4)}(H_i)]=  \mathcal{B}(\mathbb{Z}_4)[\{\mathbb{Z}_4\} \sqcup \Nbd_{\mathcal{B}(\mathbb{Z}_4)}(\mathbb{Z}_4)] = K_{1, 12}$. For  $17 \leq i \leq 20$, we have $H_i \cong \mathbb{Z}_2 \times \mathbb{Z}_2$ and so $\mathcal{B}(S_4)[\{H_i\} \sqcup \Nbd_{\mathcal{B}(S_4)}(H_i)]= \mathcal{B}(\mathbb{Z}_2 \times \mathbb{Z}_2)[\{\mathbb{Z}_2 \times \mathbb{Z}_2\} \sqcup \Nbd_{\mathcal{B}(\mathbb{Z}_2 \times \mathbb{Z}_2)}(\mathbb{Z}_2 \times \mathbb{Z}_2)] = K_{1, 6}$.
	For $21 \leq i \leq 24$, we have $H_i \cong S_3$ and so $\mathcal{B}(S_4)[\{H_i\} \sqcup \Nbd_{\mathcal{B}(S_4)}(H_i)]= \mathcal{B}(S_3)[\{S_3\} \sqcup \Nbd_{\mathcal{B}(S_3)}(S_3)] = K_{1, 18}$.
	For $25 \leq i \leq 27$, we have $H_i \cong D_8$ and so $\mathcal{B}(S_4)[\{H_i\} \sqcup \Nbd_{\mathcal{B}(S_4)}(H_i)]= \mathcal{B}(D_8)[\{D_8\} \sqcup \Nbd_{\mathcal{B}(D_8)}(D_8)] =  K_{1, 24}$.
	We have $H_{28} \cong A_4$ and so $\mathcal{B}(S_4)[\{H_{28}\} \sqcup \Nbd_{\mathcal{B}(S_4)}(H_{28})]= \mathcal{B}(A_4)[\{A_4\} \sqcup \Nbd_{\mathcal{B}(A_4)}(A_4)]= K_{1, 96}$.
	Lastly, for the subgroup $S_4$ we have $\Nbd_{\mathcal{B}(S_4)}(S_4) = S_4 \times S_4 \setminus (\sqcup_{i=0}^{28}\Nbd_{\mathcal{B}(S_4)}(H_i))$ and $\mathcal{B}(S_4)[\{S_4\} \sqcup \Nbd_{\mathcal{B}(S_4)}(S_4)] = K_{1, 216}$ noting that $|\Nbd_{\mathcal{B}(S_4)}(S_4)| = 576 - 360 = 216$.  Hence, 
	\begin{center}
		$\mathcal{B}(S_4) = K_2 \sqcup 9K_{1, 3} \sqcup 4K_{1, 8} \sqcup 3K_{1, 12} \sqcup 4K_{1, 6} \sqcup 4K_{1, 18} \sqcup 3K_{1, 24} \sqcup K_{1, 96} \sqcup K_{1, 216}$.
	\end{center}
	
	\subsection{Realization of $\mathcal{B}(D_{2p})$ and $\mathcal{B}(D_{2p^2})$}  
	In this section, we realize the graph structures of $\mathcal{B}(G)$ for the dihedral groups $D_{2p}$ and $D_{2p^2}$ where $p$ is a prime number. 
	Let us  begin with the group $D_{2p}$.
	\begin{theorem}\label{structure_of_D_2p}
		Let $D_{2p}=\langle a, b: a^p=b^2=1, bab=a^{-1} \rangle$ be the dihedral group of order $2p$, where $p$ is a prime. Then 
		$
		\mathcal{B}(D_{2p})=K_2 \sqcup pK_{1, 3} \sqcup K_{1, p^2-1} \sqcup K_{1, 3p(p-1)}.
		$
	\end{theorem}
	\begin{proof}
		\textbf{Case 1.}  $p=2$.
		
		We have $D_4= \langle a, b: a^2=b^2=1, bab=a^{-1} \rangle = \{1, a, b, ab \}$. The subgroups of $D_4$ are $H_1=\{1\}, H_2=\{1, a\}, H_3=\{1, b\}, H_4=\{1, ab\}$ and $H_5=D_4$.  
		Clearly, $(1, 1)$ is the only vertex   adjacent to  $H_1$ and so $\deg_{\mathcal{B}(D_4)}(H_1)=1$. Therefore, the  subgraph induced by  $\Nbd_{\mathcal{B}(D_4)}(H_1) \sqcup \{H_1\}$ in $\mathcal{B}(D_4)$  is $K_{1, 1}=K_2$.
		The vertices  adjacent to $H_2$ are  $(1, a), (a, 1)$ and $(a, a)$. Therefore,  $\deg_{\mathcal{B}(D_4)}(H_2)=3$ and so the  subgraph induced by  $\Nbd_{\mathcal{B}(D_4)}(H_2) \sqcup \{H_2\}$ in $\mathcal{B}(D_4)$  is $K_{1, 3}$. Similarly, $H_{3}$ is adjacent to $(1, b), (b, 1), (b, b)$ and $H_4$ is adjacent to $(1, ab), (ab, 1), (ab, ab)$. So, $\deg_{\mathcal{B}(D_4)}(H_3)=\deg_{\mathcal{B}(D_4)}(H_4)=3$ and 
		\[
		\mathcal{B}(D_4)[\Nbd_{\mathcal{B}(D_4)}(H_3) \sqcup \{H_3\}] = K_{1, 3} =   \mathcal{B}(D_4)[\Nbd_{\mathcal{B}(D_4)}(H_4) \sqcup \{H_4\}].
		\]
		
		Lastly, $H_5$ is adjacent to \quad $(a, b), (b, a), (a, ab), (ab, a), (b, ab)$ and $(ab , b)$. Therefore,  $\deg_{\mathcal{B}(D_4)}(H_5)$ $=6$ and so the subgraph induced by $\Nbd_{\mathcal{B}(D_4)}(H_5) \sqcup \{H_5\}$ in $\mathcal{B}(D_4)$  is $K_{1, 6}$. Thus,
		\begin{align*}
			\mathcal{B}(D_4) &= \underset{H \in L(D_4)}{\sqcup}\mathcal{B}(D_4)[\Nbd_{\mathcal{B}(D_4)}(H) \sqcup \{H\}]\\
			&= K_2 \sqcup 3K_{1, 3} \sqcup K_{1, 6} 
			= K_2 \sqcup 2K_{1, 3} \sqcup K_{1, 2^2-1} \sqcup K_{1, 3\cdot2(2-1)}.
		\end{align*}
		\textbf{Case 2.}  $p$ is an odd prime.
		
		The subgroups of $D_{2p}=\langle a, b: a^p=b^2=1, bab=a^{-1} \rangle=\{1, a, a^2, a^3, \ldots, a^{p-1}, b, ab,$ $a^2b, \ldots, a^{p-1}b\}$ are $H_0=\{1\}, H_1=\{1, b\}, H_2=\{1, ab\}, H_3=\{1, a^2b\}, H_4=\{1, a^3b\},$ $\ldots, H_p=\{1, a^{p-1}b\}, H_{p+1}=\{1, a, a^2, \ldots, a^{p-1}\}$ and $H_{p+2}=D_{2p}$.  
		Clearly, $(1, 1)$  is the only vertex adjacent to $H_0$ and  so $\deg_{\mathcal{B}(D_{2p})}(H_0)=1$. Therefore, the subgraph induced by  $\Nbd_{\mathcal{B}(D_{2p})}(H_0) \sqcup \{H_0\}$ in $\mathcal{B}(D_{2p})$  is $K_2$. The vertices  adjacent to $H_1$ are  $(1, b), (b, 1), (b, b)$. Therefore,  $\deg_{\mathcal{B}(D_{2p})}(H_1)=3$ and so the subgraph induced by  $\Nbd_{\mathcal{B}(D_{2p})}(H_1) \sqcup \{H_1\}$ in $\mathcal{B}(D_{2p})$ is $K_{1, 3}$. Similarly, for each $i = 2, 3, \dots, p$ the  subgraph induced by $\Nbd_{\mathcal{B}(D_{2p})}(H_i) \sqcup \{H_i\}$ in $\mathcal{B}(D_{2p})$  is $K_{1, 3}$. 
		The vertices adjacent to $H_{p+1}$ are $(1, a), (1, a^2), \ldots, (1, a^{p-1}), (a, 1),$ $ (a^2, 1), \ldots,$ $(a^{p-1}, 1)$ and $(a^i, a^j)$ where $1 \leq i, j \leq p-1$.
		Therefore, $\deg_{\mathcal{B}(D_{2p})}(H_{p+1})=(p-1)+(p-1)+(p-1)^2=(p-1)(p+1)=p^2-1$ and so the  subgraph induced by $\Nbd_{\mathcal{B}(D_{2p})}(H_{p+1}) \sqcup \{H_{p+1}\}$ in $\mathcal{B}(D_{2p})$  is $K_{1, p^2-1}$. Lastly, the vertices adjacent to $H_{p+2}$ are $(a^i, a^jb)$, where $1 \leq i \leq p-1$ and $0 \leq j \leq p-1$; $(a^jb, a^i)$, where $1 \leq i \leq p-1$ and $0 \leq j \leq p-1$ and $(a^ib, a^jb)$, where $0 \leq i\ne j \leq p-1$. 
		Therefore, $\deg_{\mathcal{B}(D_{2p})}(H_{p+2})=p(p-1)+p(p-1)+p^2-p=3p(p-1)$ and the  subgraph induced by  $\Nbd_{\mathcal{B}(D_{2p})}(H_{p+2}) \sqcup \{H_{p+2}\}$ in $\mathcal{B}(D_{2p})$ is $K_{1, 3p(p-1)}$. Thus, 
		\begin{align*}
			\mathcal{B}(D_{2p}) &= \underset{H \in L(D_{2p})}{\sqcup}\mathcal{B}(D_{2p})[\Nbd_{\mathcal{B}(D_{2p})}(H) \sqcup \{H\}]\\
			&= K_2 \sqcup pK_{1, 3} \sqcup K_{1, p^2-1} \sqcup K_{1, 3p(p-1)}.
		\end{align*}
		This completes the proof.
	\end{proof}
	\begin{theorem}\label{structure_of_D_2p2}
		Let $D_{2p^2}=\langle a, b: a^{p^2}=b^2=1, bab=a^{-1} \rangle$ be the dihedral group of order $2p^2$, where $p$ is a prime. Then 
		\[
		\mathcal{B}(D_{2p^2})=K_2 \sqcup p^2K_{1, 3} \sqcup K_{1, p^2-1} \sqcup K_{1, p^4-p^2} \sqcup pK_{1, 3p(p-1)} \sqcup K_{1, 3p^2(p^2-p)}.
		\]
	\end{theorem}
	\begin{proof}
		If $p = 2$ then we have already obtained that $\mathcal{B}(D_8)= K_2 \sqcup 5K_{1, 3} \sqcup 2K_{1, 6} \sqcup K_{1, 12} \sqcup K_{1, 24}$. Therefore, we consider the case when $p$ is an odd prime.
		The subgroups of $D_{2p^2}= \langle a, b: a^{p^2}=b^2=1, bab=a^{-1} \rangle = \{1, a, a^2, \ldots a^{p^2-1}, b, ab, a^2b,$  $\ldots, a^{p^2-1}b\}$ are
		
		\begin{center}
			$I=\{1\}$; \quad
			$H_i=\langle a^ib \rangle=\{1, a^ib\}$ for  $0 \leq i \leq p^2-1$; 
			
			$K=\langle a^p \rangle=\{1, a^p, a^{2p}, \ldots, a^{(p-1)p}\}$; \quad
			$T=\langle a \rangle=\{1, a, a^2, \ldots, a^{p^2-1}\}$;

			\noindent $M_r=\langle a^p, a^rb \rangle=\{1, a^p, a^{2p}, \ldots, a^{(p-1)p}, a^{p+r}b, a^{2p+r}b, a^{3p+r}b, \ldots, a^{(p-1)p+r}b, a^rb\}$ for $0 \leq r \leq p-1$; 
			and $G=D_{2p^2}$. 
		\end{center}
		
		Thus, $L(D_{2p^2})=\{I, H_0, H_1, \ldots, H_{p^2-1}, K, T, M_0, M_1, M_2, \ldots, M_{p-1}, G\}$. 
		Clearly, $(1, 1)$ is the   only vertex adjacent to $I$ and  so $\deg_{\mathcal{B}(D_{2p^2})}(I)=1$. Therefore, $\mathcal{B}(D_{2p^2})[\Nbd_{\mathcal{B}(D_{2p^2})}(I)$ $\sqcup \{I\}] = K_2$. We have
		$\Nbd_{\mathcal{B}(D_{2p^2})}(H_i) = \{(1, a^ib), (a^ib, 1), (a^ib, a^ib)\}$ and so $\deg_{\mathcal{B}(D_{2p^2})}(H_i)$ $=3$ for $0 \leq i \leq p^2-1$. Therefore, $\mathcal{B}(D_{2p^2})[\Nbd_{\mathcal{B}(D_{2p^2})}(H_i) \sqcup \{H_i\}] = K_{1, 3}$ for $0 \leq i \leq p^2-1$.  We have 
		$\Nbd_{\mathcal{B}(D_{2p^2})}(K) =\{(a^{ip}, a^{jp}): 0 \leq i, j \leq p-1\} \setminus \{(1, 1)\}$ and so $\deg_{\mathcal{B}(D_{2p^2})}(K)=p^2-1$. Therefore, $\mathcal{B}(D_{2p^2})[\Nbd_{\mathcal{B}(D_{2p^2})}(K)\sqcup \{K\}] = K_{1, p^2-1}$. 
		For the subgroup $T$, we have  
		$\Nbd_{\mathcal{B}(D_{2p^2})}(T) = \{(a^r, a^s): \langle a^r, a^s \rangle = T\}$. 
		We know that $\langle a^r, a^s \rangle = \langle a^{\gcd(r,s)} \rangle$ and $\langle a^l \rangle =T$ if and only if $\gcd(l, p^2)=1$. Therefore, $|\Nbd_{\mathcal{B}(D_{2p^2})}(T)| = p^4-p^2 = \deg_{\mathcal{B}(D_{2p^2})}(T)$
		and so $\mathcal{B}(D_{2p^2})[\Nbd_{\mathcal{B}(D_{2p^2})}(T)\sqcup \{T\}] = K_{1, p^4-p^2}$. 
		Finally, for the subgroups $M_r$, we have
		$\Nbd_{\mathcal{B}(D_{2p^2})}(M_r) = \{(a^{ip}b, a^{jp+r}), (a^{jp+r}, a^{ip}b), (a^{ip+r}b, a^{jp+r}) : 0 \leq i \ne j \leq p-1\}$ for $0 \leq r \leq p-1$.
		Therefore, $|\Nbd_{\mathcal{B}(D_{2p^2})}(M_r)| = 2p(p-1)+p^2-p=3p(p-1) = \deg_{\mathcal{B}(D_{2p^2})}(M_r)$ and so $\mathcal{B}(D_{2p^2})[\Nbd_{\mathcal{B}(D_{2p^2})}(M_r)\sqcup \{M_r\}] = K_{1, 3p(p-1)}$ for $0 \leq r \leq p-1$.
		
		By Lemma \ref{deg_sum=num_of_edges}, we get $\deg_{\mathcal{B}(D_{2p^2})}(D_{2p^2}) = 4p^4-(1+3p^2+p^2-1+p^4-p^2+3p^3-3p^2)=3p^2(p^2-p)$. Therefore, $\mathcal{B}(D_{2p^2})[\Nbd_{\mathcal{B}(D_{2p^2})}(G) \sqcup \{G\}] = K_{1, 3p^2(p^2-p)}$. Hence, 
		\begin{align*}
			\mathcal{B}(D_{2p^2})  &= \underset{H \in L(D_{2p^2})}{\sqcup}\mathcal{B}(D_{2p^2})[\Nbd_{\mathcal{B}(D_{2p^2})}(H) \sqcup \{H\}]\\	
			&=K_2 \sqcup p^2K_{1, 3} \sqcup K_{1, p^2-1} \sqcup K_{1, p^4-p^2} \sqcup pK_{1, 3p(p-1)} \sqcup K_{1, 3p^2(p^2-p)}. 
		\end{align*}
		This completes the proof.
	\end{proof}
	
	\subsection{Realization of $\mathcal{B}(Q_{4p})$ and $\mathcal{B}(Q_{4p^2})$} 
	In this section, we realize the graph structures of $\mathcal{B}(G)$ for the dicyclic groups $Q_{4p}$ and $Q_{4p^2}$, where $p$ is a prime number.
	We  begin with the group $Q_{4p}$.
	\begin{theorem}
		Let $Q_{4p} =  \langle a, b : a^{2p} = 1, b^2 = a^p, bab^{-1} = a^{-1} \rangle$ be the dicyclic group of order $4p$, where $p$ is a prime. Then
		\[
		\mathcal{B}(Q_{4p})=\begin{cases}
			K_2 \sqcup K_{1, 3} \sqcup 3K_{1, 12} \sqcup K_{1, 24}, & \text{ when } p=2 \\
			K_2 \sqcup K_{1, 3} \sqcup pK_{1, 12} \sqcup K_{1, p^2-1} \sqcup K_{1, 3p^2-3} \sqcup K_{1, 12p^2-12p}, & \text{ when } p \geq 3.
		\end{cases}
		\]
	\end{theorem}
	\begin{proof}
		If $p=2$ then we have already obtained that $\mathcal{B}(Q_8)=K_2 \sqcup K_{1, 3} \sqcup 3K_{1, 12} \sqcup K_{1, 24}$. Therefore, we consider the case when  $p$ is an odd prime.
		
		The subgroups of $Q_{4p}=  \langle a, b : a^{2p} = 1, b^2 = a^p, bab^{-1} = a^{-1} \rangle$ are $I=\{1\}$, $K=\{1, a^p\}$, $T=\langle a^2 \rangle =\{a^2, a^4, \ldots, a^{2p}=1\}$, $S=\langle a \rangle =\{a, a^2, \ldots, a^{2p}=1\}$, $H_i=\langle a^ib \rangle=\{1, a^ib, b^2, a^{p+i}b\}$ for $1 \leq i \leq p$; and $G=Q_{4p}$. Thus, $L(Q_{4p})=\{I, K, T, S, H_1, H_2, \ldots, H_p,$ $ G\}$. Clearly, $(1, 1)$ is the only vertex adjacent to $I$. Therefore, $\mathcal{B}(Q_{4p})[\{I\} \sqcup \Nbd_{\mathcal{B}(Q_{4p})}(I)]$ $=K_2$. Since $K, T$ and $S$ are cyclic subgroups of order two, $p$ and $2p$ respectively,  by Observation \ref{vrtex_deg_of_X_in_B(G)}(b), we have $\mathcal{B}(Q_{4p})[\{K\} \sqcup \Nbd_{\mathcal{B}(Q_{4p})}(K)]=K_{1, 3}$, $\mathcal{B}(Q_{4p})[\{T\} \sqcup \Nbd_{\mathcal{B}(Q_{4p})}(T)]=K_{1, p^2-1}$ and $\mathcal{B}(Q_{4p})[\{S\} \sqcup \Nbd_{\mathcal{B}(Q_{4p})}(S)]=K_{1, 3p^2-3}$. Also, since $H_i$'s are cyclic subgroups of order four,  by Observation \ref{vrtex_deg_of_X_in_B(G)}(b), we have $\mathcal{B}(Q_{4p})[\{H_i\} \sqcup \Nbd_{\mathcal{B}(Q_{4p})}(H_i)]=K_{1, 12}$ for $1 \leq i \leq p$. 
		
		By Lemma \ref{deg_sum=num_of_edges}, we get $\deg_{\mathcal{B}(Q_{4p})}(Q_{4p})=|\Nbd_{\mathcal{B}(Q_{4p})}(Q_{4p})|=16p^2-(1+3+p^2-1+3p^2-3+12p)=12p^2-12p$. Therefore, $\mathcal{B}(Q_{4p})[\{G\} \sqcup \Nbd_{\mathcal{B}(Q_{4p})}(G)]=K_{1, 12p^2-12p}$. Hence,
		\begin{align*}
			\mathcal{B}(Q_{4p})&=\underset{H \in L(Q_{4p})}{\sqcup} \mathcal{B}(Q_{4p})[\{H\} \sqcup \Nbd_{\mathcal{B}(Q_{4p})}(H)] \\
			&=K_2 \sqcup K_{1, 3} \sqcup pK_{1, 12} \sqcup K_{1, p^2-1} \sqcup K_{1, 3p^2-3} \sqcup K_{1, 12p^2-12p}.
		\end{align*}
		This completes the proof.
	\end{proof}
	\begin{theorem}
		Let $Q_{4p^2} =  \langle a, b : a^{2p^2} = 1, b^2 = a^{p^2}, bab^{-1} = a^{-1} \rangle$ be the dicyclic group of order $4p^2$, where $p$ is a prime. Then
		\[
		\mathcal{B}(Q_{4p^2})=\begin{cases}
			K_2 \sqcup K_{1, 3} \sqcup 5K_{1, 12} \sqcup 2K_{1, 24} \sqcup K_{1, 48} \sqcup K_{1, 96}, & \text{ when } p=2 \\
			K_2 \sqcup K_{1, 3} \sqcup p^2K_{1, 12} \sqcup K_{1, p^2-1} \sqcup K_{1, 3p^2-3} \sqcup K_{1, 3p^4-3p^2} \\ \qquad \qquad \sqcup (p-1)K_{1, 12p^2-12p} \sqcup K_{1, 13p^4-12p^3+11p^2-12p}, & \text{ when } p \geq 3.
		\end{cases}
		\]
	\end{theorem}
	\begin{proof}
		\textbf{Case 1.} $p=2$.
		
		We have $Q_{16}=\{1, a, a^2, \ldots, a^7, b, ab, a^2b, \ldots, a^7b\}$. The subgroups of $Q_{16}$ are 
		\begin{center}
			$I=\{1\}$; \quad $J=\langle a^4 \rangle =\{1, a^4\}$; \quad $K=\langle a^2 \rangle=\{1, a^2, a^4, a^6\}$; $L=\langle a \rangle=\{1, a, a^2, \ldots, a^7\}$; \quad $H_i=\langle a^ib \rangle =\{1, a^ib, a^4, a^{p^2+i}b\}$ for $1 \leq i \leq 4$; $M_x=\langle a^2, x: (a^2)^2=x^2, (a^2)^4=1, xa^2x^{-1}=(a^2)^{-1} \rangle$ for $x=b$ and $ab$; 
			
			and $G=Q_{16}$.
		\end{center} 
		Thus, $L(Q_{16})=\{I, J, K, L, H_1, \ldots, H_4, M_b, M_{ab}, G\}$. Clearly, $(1, 1)$ is the only vertex adjacent to $I$. Therefore, $\mathcal{B}(Q_{16})[\{I\} \sqcup \Nbd_{\mathcal{B}(Q_{16})}(I)]=K_2$. Since $J \cong \mathbb{Z}_2$, by Observation \ref{vrtex_deg_of_X_in_B(G)}(b), $\mathcal{B}(Q_{16})[\{J\} \sqcup \Nbd_{\mathcal{B}(Q_{16})}(J)]=K_{1, 3}$. Also, for $1 \leq i \leq 4$, $K \cong H_i \cong \mathbb{Z}_4$ and so $\mathcal{B}(Q_{16})[\{K\} \sqcup \Nbd_{\mathcal{B}(Q_{16})}(K)]=\mathcal{B}(Q_{16})[\{H_i\} \sqcup \Nbd_{\mathcal{B}(Q_{16})}(H_i)]=K_{1, 12}$. Now, $\Nbd_{\mathcal{B}(Q_{16})}(L)=L \times L \setminus (\Nbd_{\mathcal{B}(Q_{16})}(I)\sqcup \Nbd_{\mathcal{B}(Q_{16})}(J)\sqcup \Nbd_{\mathcal{B}(Q_{16})}(K))$. That is, $|\Nbd_{\mathcal{B}(Q_{16})}(L)|=64-16=48$ and so $\mathcal{B}(Q_{16})[\{L\} \sqcup \Nbd_{\mathcal{B}(Q_{16})}(L)]=K_{1, 48}$. Also, $M_b \cong M_{ab} \cong Q_8$ and so $\mathcal{B}(Q_{16})[\{M_b\} \sqcup \{M_{ab}\} \sqcup \Nbd_{\mathcal{B}(Q_{16})}(M_a) \sqcup \Nbd_{\mathcal{B}(Q_{16})}(M_{ab})]=2K_{1, 24}$.
		
		By Lemma \ref{deg_sum=num_of_edges}, we get $\deg_{\mathcal{B}(Q_{16})}(Q_{16})=|\Nbd_{\mathcal{B}(Q_{16})}(Q_{16})|=256-(1+3+60+48+48)=96$. Therefore, $\mathcal{B}(Q_{16})[\{G\} \sqcup \Nbd_{\mathcal{B}(Q_{16})}(G)]=K_{1, 96}$. Hence,
		\begin{align*}
			\mathcal{B}(Q_{16})&=\underset{H \in L(Q_{16})}{\sqcup} \mathcal{B}(Q_{16})[\{H\} \sqcup \Nbd_{\mathcal{B}(Q_{16})}(H)] \\
			&=K_2 \sqcup K_{1, 3} \sqcup 5K_{1, 12} \sqcup 2K_{1, 24} \sqcup K_{1, 48} \sqcup K_{1, 96}.
		\end{align*}
		\textbf{Case 2.} $p$ is an odd prime.
		
		The subgroups of $Q_{4p^2}=\{1, a, a^2, \ldots, a^{2p^2-1}, b, ab, a^2b, \ldots, a^{2p^2-1}b\}$ are
		\begin{center}
			$I=\{1\}$; \quad $J=\{1, a^{p^2}\}$; \quad $K=\langle a^{2p} \rangle=\{a^{2p}, a^{4p}, \ldots, (a^{2p})^p=1\}$; $L=\langle a^p \rangle=\{a^p, a^{2p}, \ldots, (a^p)^{2p}=1\}$; \quad $T=\langle a \rangle=\{a, a^2, \ldots, a^{2p^2}=1\}$; $H_i=\langle a^ib \rangle=\{a^ib, a^{p^2}, a^{p^2+i}b, 1\}$ for $1 \leq i \leq p^2$; $M_{x_i}=\langle (a^p)^i, x: ((a^p)^i)^p=x^2, ((a^p)^i)^{2p}=1, x(a^p)^ix^{-1}=((a^p)^i)^{-1}\rangle$, where $x=b$ and $ab$, $i < p$ and $\gcd(i, 2p)=1$; and $G=Q_{4p^2}$.
		\end{center}
		Thus, \quad $L(Q_{4p^2}) \quad =\{I, J, K, L, T, H_1, H_2, \ldots, H_{p^2}, M_{b_1}, M_{b_3} \ldots, M_{b_{p-2}}, M_{ab_1}, M_{ab_3}, \ldots,$ 
		
		\noindent $ M_{ab_{p-2}}, G\}$. Clearly, $(1, 1)$ is the only vertex adjacent to $I$. Therefore, $\mathcal{B}(Q_{4p^2})[\{I\} \sqcup \Nbd_{\mathcal{B}(Q_{4p^2})}(I)]=K_2$. Since $J \cong \mathbb{Z}_2$, $K \cong \mathbb{Z}_p$, $L \cong \mathbb{Z}_{2p}$, $T \cong \mathbb{Z}_{2p^2}$ and $H_i \cong \mathbb{Z}_4$ for $1 \leq i \leq p^2$,  by Observation \ref{vrtex_deg_of_X_in_B(G)}(b), we have
		\begin{center}
			$\mathcal{B}(Q_{4p^2})[\{J\} \sqcup \Nbd_{\mathcal{B}(Q_{4p^2})}(J)]=K_{1, 3}$; \quad $\mathcal{B}(Q_{4p^2})[\{K\} \sqcup \Nbd_{\mathcal{B}(Q_{4p^2})}(K)]=K_{1, p^2-1}$; \quad $\mathcal{B}(Q_{4p^2})[\{L\} \sqcup \Nbd_{\mathcal{B}(Q_{4p^2})}(L)]=K_{1, 3p^2-3}$; \quad $\mathcal{B}(Q_{4p^2})[\{T\} \sqcup \Nbd_{\mathcal{B}(Q_{4p^2})}(T)]=K_{1, 3p^4-3p^2}$; \quad and $\mathcal{B}(Q_{4p^2})[\{H_i\} \sqcup \Nbd_{\mathcal{B}(Q_{4p^2})}(H_i)]=K_{1, 12}$ for $1 \leq i \leq p^2$.
		\end{center}
		Also, for odd $i$, $1 \leq i \leq p-2$, we have $M_{b_i} \cong M_{ab_i} \cong Q_{4p}$. As such, $\mathcal{B}(Q_{4p^2})[\{M_{b_i}\} \sqcup \Nbd_{\mathcal{B}(Q_{4p^2})}(M_{b_i})]=K_{1, 12p^2-12p}=\mathcal{B}(Q_{4p^2})[\{M_{ab_i}\} \sqcup \Nbd_{\mathcal{B}(Q_{4p^2})}(M_{ab_i})]$ for odd $i$, $1 \leq i \leq p-2$. By Lemma \ref{deg_sum=num_of_edges}, we get $\deg_{\mathcal{B}(Q_{4p^2})}(Q_{4p^2})=|\Nbd_{\mathcal{B}(Q_{4p^2})}(Q_{4p^2})|=16p^4-(1+3+p^2-1+3p^2-3+3p^4-3p^2+12p^2+12p^3-12p^2-12p^2+12p)=13p^4-12p^3+11p^2-12p$. Therefore, $\mathcal{B}(Q_{4p^2})[\{G\} \sqcup \Nbd_{\mathcal{B}(Q_{4p^2})}(G)]=K_{1, 13p^4-12p^3+11p^2-12p}$. Hence, 
		\begin{align*}
			\mathcal{B}(Q_{4p^2})&=\underset{H \in L(Q_{4p^2})}{\sqcup} \mathcal{B}(Q_{4p^2})[\{H\} \sqcup \Nbd_{\mathcal{B}(Q_{4p^2})}(H)] \\
			&=K_2 \sqcup K_{1, 3} \sqcup p^2K_{1, 12} \sqcup K_{1, p^2-1} \sqcup K_{1, 3p^2-3} \sqcup K_{1, 3p^4-3p^2} \\
			& \qquad \qquad \qquad \sqcup (p-1)K_{1, 12p^2-12p} \sqcup K_{1, 13p^4-12p^3+11p^2-12p}.
		\end{align*}
		This completes the proof.
	\end{proof}

	It may be interesting to realize the structures of $\mathcal{B}(G)$ for the dihedral group and dicyclic group (in general), 
	quasi-dihedral group $QD_{2^{n}} = \langle a, b : a^{2^n} = b^2 = 1, bab^{-1} = a^{-1} \rangle$ ($n\geq 3$) etc. We leave this as a problem for further research. In a separate paper we shall compute various spectra, energies and certain topological indices of $\mathcal{B}(G)$  for the groups considered in this paper.

	\section*{Acknowledgements}
	The authors would like to thank the referee for his/her valuable comments. The first author is thankful to Council of Scientific and Industrial Research  for the fellowship (File No. 09/0796(16521)/2023-EMR-I).

\end{document}